\documentclass[12pt,reqno]{amsart}

\usepackage[T1]{fontenc}
\usepackage{lmodern}

\usepackage{amsmath,amssymb,amsthm,mathrsfs,comment}
\usepackage{amsfonts} 
\usepackage{mathtools}
\usepackage{enumitem}
\usepackage{xspace}
\usepackage{bbold}
\usepackage{amscd}
\usepackage{verbatim}

\usepackage{marvosym}

\usepackage{xcolor}
\usepackage{bbold}

\usepackage{mathpazo}
\usepackage[scaled=0.95]{helvet}
\usepackage{courier}

\usepackage[
  a4paper,
  left=3cm,
  right=2.5cm,
  top=2.8cm,
  bottom=3.8cm,
  headheight=15pt,
  headsep=0.8cm,
  footskip=1.5cm
]{geometry}

\usepackage[
    colorlinks,
    allcolors=black,
    pdfauthor={},
    pdftitle={}
]{hyperref}
\usepackage[nameinlink]{cleveref}

\usepackage{enumitem}

\setlist[enumerate,1]{
  label={\upshape(\roman*)},
}

\setlist[itemize,1]{
  leftmargin=*,
  topsep=6pt,
  itemsep=2pt
}

\theoremstyle{plain}
\newtheorem{theorem}{Theorem}[section]
\newtheorem{proposition}[theorem]{Proposition}
\newtheorem{corollary}[theorem]{Corollary}
\newtheorem{lemma}[theorem]{Lemma}

\theoremstyle{definition}
\newtheorem{definition}[theorem]{Definition}
\newtheorem{example}[theorem]{Example}

\DeclareSymbolFont{bbold}{U}{bbold}{m}{n}
\DeclareSymbolFontAlphabet{\mathbbold}{bbold}
\newcommand{\one}{\mathbbold{1}}

\DeclareMathOperator{\Span}{span}

\DeclareMathOperator{\ac}{AC}
\DeclareMathOperator{\bv}{BV}
\DeclareMathOperator{\cbv}{CBV}

\DeclareMathOperator{\Log}{Log}

\DeclareMathOperator{\aut}{Aut}
\DeclareMathOperator{\isom}{Isom}

\newcommand{\Z}{\mathbb{Z}}
\newcommand{\R}{\mathbb{R}}

\usepackage[
  backend=biber,
  style=alphabetic,
  sorting=nyt,
  maxbibnames=99,
  giveninits=true,
  doi=false,
  url=false,
  isbn=false,
  eprint=true,
  backref=false
]{biblatex}

\DeclareFieldFormat{pages}{#1}
\DeclareFieldFormat[article]{pages}{#1}
\DeclareFieldFormat[incollection]{pages}{#1}
\DeclareFieldFormat[inproceedings]{pages}{#1}

\DeclareFieldFormat[article]{journaltitle}{\mkbibemph{#1}}
\DeclareFieldFormat[book]{title}{\mkbibemph{#1}}
\DeclareFieldFormat[incollection]{title}{#1}
\DeclareFieldFormat[phdthesis]{title}{\mkbibemph{#1}}

\newcommand*\mkparenyear{%
  \iffieldundef{year}
    {}
    {\mkbibparens{\printfield{year}}}%
}

\DeclareBibliographyDriver{article}{%
  \printnames{author}%
  \setunit{\addcomma\space}%
  \printfield[titlecase]{title}%
  \setunit{\addcomma\space}%
  \printfield{journaltitle}%
  \setunit{\addspace}%
  \printfield{volume}%
\space
\mkparenyear%
  \setunit{\addcomma\space}%
  \printfield{pages}%
  \finentry
}
\DeclareBibliographyDriver{book}{%
  \ifnameundef{author}
    {\printnames{editor}\setunit{\addspace}\printtext[parens]{\bibstring{editors}}}
    {\printnames{author}}%
  \setunit{\addcomma\space}%
  \printfield{title}%
  \setunit{\addcomma\space}%
  \printfield{edition}%
  \setunit{\addcomma\space}%
  \printlist{publisher}%
  \setunit{\addcomma\space}%
  \printlist{location}%
  \setunit{\addcomma\space}%
  \printfield{year}%
  \finentry
} 
\DeclareBibliographyDriver{incollection}{%
  \printnames{author}%
  \setunit{\addcomma\space}%
  \printfield{title}%
  \setunit{\addcomma\space}%
  \printfield{booktitle}%
  \setunit{\addcomma\space}%
  \printlist{publisher}%
  \setunit{\addcomma\space}%
  \printlist{location}%
  \setunit{\addcomma\space}%
  \printfield{year}%
  \finentry
}

\DeclareBibliographyDriver{misc}{%
  \printnames{author}%
  \setunit{\addcomma\space}%
  \printfield{title}%
  \setunit{\addcomma\space}%
  \printfield{note}%
  \setunit{\addcomma\space}%
  \printfield{year}%
  \finentry
}

\DeclareBibliographyDriver{thesis}{%
  \printnames{author}%
  \setunit{\addcomma\space}%
  \printfield[noformat]{title}%
  \setunit{\addcomma\space}%
  \printfield{type}%
  \setunit{\addcomma\space}%
  \printlist{institution}%
  \setunit{\addcomma\space}%
  \printfield{year}%
  \finentry
}

\DeclareBibliographyDriver{phdthesis}{%
  \printnames{author}%
  \setunit{\addcomma\space}%
  \printfield{title}%
  \setunit{\addcomma\space}%
  \printlist{school}%
  \setunit{\addcomma\space}%
  \printfield{year}%
  \finentry
}

\begin{document}

\title{Isometries for Hilbert's metric on vector lattices of continuous functions}
 
\author{Marko Kandi\'c}
\address{Faculty of Mathematics and Physics, University of Ljubljana,
  Jadranska 19, 1000 Ljubljana, Slovenia \ \ and \ \ \
  Institute of Mathematics, Physics, and Mechanics, Jadranska 19, 1000 Ljubljana, Slovenia}
\email{marko.kandic@fmf.uni-lj.si}

\author{Mark Roelands}
\address{Mathematical Institute, Leiden University,
  Einsteinweg 55,
2333 CC, Leiden, the Netherlands}
\email{m.roelands@math.leidenuniv.nl}



\date{\today}
 
\begin{abstract}
In this paper Hilbert's projective metric on the interior of the positive cone of an Archimedean vector lattice with a strong order unit is studied. We show that every surjective isometry between two such projective cones is induced by a lattice isomorphism, possibly composed with the inversion map, and that the inversion occurs precisely when the vector lattice is a bounded inverse closed $f$-algebra. In particular, two such projective cones are isometric if and only if the underlying vector lattices are lattice isomorphic, and we determine the corresponding isometry groups. We apply these results to spaces of H\"older continuous functions, absolutely continuous functions, and continuous functions of bounded variation, for which we describe the surjective isometries explicitly. Finally, we show that the projective cones of these spaces of functions on an interval all have the same metric completion, but that no two of them are isometric.
\end{abstract}
 
\maketitle

\section{Introduction}

In \cite{Hil95}, Hilbert introduced a metric on a bounded open convex set $\Omega$ in a finite dimensional real vector space. For distinct $x,y\in\Omega$, let $x'$ and $y'$ be the points where the straight line through $x$ and $y$ meets the boundary of $\Omega$, labelled such that $x$ lies between $x'$ and $y$, and $y$ lies between $x$ and $y'$. Hilbert's metric is then given by
\[
\delta_H(x,y):=\log\frac{\|x'-y\|\,\|y'-x\|}{\|x'-x\|\,\|y'-y\|}.
\]
The resulting metric spaces, which generalize Klein's model of hyperbolic space, are known as Hilbert geometries; see \cite{PT14} for an overview. In his work on extensions of Jentzsch's theorem, Birkhoff \cite{Bir57} introduced a version of Hilbert's metric on the rays in the interior of a cone, which is defined in terms of the partial order induced by the cone. If the cone admits a strictly positive linear functional, then Birkhoff's version coincides with Hilbert's original metric on the corresponding cross-section of the cone \cite[Theorem~2.1.2]{LN12}. Birkhoff showed that positive linear maps are contractions for this metric, and his version of Hilbert's metric has since become a basic tool in the spectral theory of linear and nonlinear operators on cones; see \cite{LN12,Nus88} and the references therein. In this paper we work with Birkhoff's version of Hilbert's metric on the interior of the positive cone of an order unit space, as recalled in Subsection~\ref{subsec:hilbert}.

The isometries of Hilbert geometries and of Hilbert's metric on cones have been studied extensively. In finite dimensions, they were determined for simplices by de la Harpe \cite{dlH93}, for polyhedral domains by Lemmens and Walsh \cite{LW11}, and in dimension two by Matveev and Troyanov \cite{MT15}; see also \cite{Spe14}. Walsh \cite{Wal18} showed that, for a finite dimensional order unit space, the group of Hilbert isometries of the interior of the cone is strictly larger than the group of projectivities if and only if the space is a Euclidean Jordan algebra that is neither $\mathbb R^2$ nor a spin factor, in which case the additional isometries arise from the inversion $x\mapsto x^{-1}$. In infinite dimensions, Moln\'ar \cite{Mol09} described the Hilbert isometries of the cone of positive invertible operators on a Hilbert space, and Bosch\'e \cite{Bos12} those of symmetric cones. Lemmens, the second author, and Wortel characterized the Hilbert isometries of infinite dimensional strictly convex domains and of the cones of  $C(K)$ in \cite{LemRoeWor18}, and of the cones in JBW-algebras \cite{LRW19}. The second author and Wortel extended the latter result to JB-algebras in \cite{RW19}. We refer to \cite{HOS84} for the theory of Jordan algebras and JB-algebras. In all these cases, the isometries are induced by linear maps, possibly composed with the inversion. Our aim is to study the isometries for Hilbert's metric on cones in Archimedean vector lattices with a strong order unit.

It was asked in \cite{LemRoeWor18} whether, in a general order unit space, the existence of a Hilbert isometry that is not induced by a linear map forces a Jordan algebra structure on the space. For complete order unit spaces, this question was answered in \cite{RoelandsTiersma2025}. Our main result, Theorem~\ref{T:d_H-isom between E and F}, describes the surjective Hilbert isometries between the projective cones $P(E_+^\circ)$ and $P(F_+^\circ)$ of Archimedean vector lattices $E$ and $F$ with strong order units, which need not be complete. Here $P(E_+^\circ)$ denotes the set of rays in the interior $E_+^\circ$ of the positive cone of $E$. These isometries are exactly the maps induced by lattice isomorphisms from $E$ onto $F$, possibly composed with the inversion, and the inversion can only occur if $E$ and $F$ are bounded inverse closed $f$-algebras in the sense of Definition~\ref{def:bic}. By Theorem~\ref{thm:f-algebra-structure-kakutani-representation}, this means that $E$ and $F$ can be represented as norm dense subalgebras of spaces of continuous functions $C(K)$ and $C(L)$ that are closed under inversion of their strictly positive elements. Since the spaces $C(K)$ are precisely the associative unital JB-algebras \cite[Theorem~3.2.2]{HOS84}, a surjective Hilbert isometry that is not induced by a lattice isomorphism exists only if the vector lattice is a dense subalgebra of an associative JB-algebra. This gives an affirmative answer to the question of \cite{LemRoeWor18} for Archimedean vector lattices with a strong order unit, where no completeness is assumed. As a consequence, the projective cones of two such vector lattices are isometric for Hilbert's metric if and only if the vector lattices are lattice isomorphic (Corollary~\ref{C: d_H isometric iff lattice isom}), which is a Banach--Stone type theorem for Hilbert's metric. We also determine the isometry group of $P(E_+^\circ)$ in Theorem~\ref{thm:isom-group-isomorphism}. It consists of the isometries induced by lattice automorphisms of $E$, together with their compositions with the inversion if $E$ is a bounded inverse closed $f$-algebra of dimension at least three.

The inversion requires a multiplication, which a vector lattice does not carry a priori. In Section~\ref{sec:bic} we show that an Archimedean vector lattice with a strong order unit admits at most one $f$-algebra multiplication with this unit, namely the one induced by any of its Kakutani representations (Theorem~\ref{thm:f-algebra-structure-kakutani-representation}). For complete spaces, this is closely related to results of Mu\~noz-Lahoz and Tradacete \cite{MunozLahozTradacete2025} and of Mu\~noz-Lahoz \cite{MunozLahoz2026} on $f$-algebra products on spaces $C(K)$.

In Section~\ref{sec:applications} we apply our results to spaces of H\"older continuous functions, absolutely continuous functions, and continuous functions of bounded variation. Up to a weight, the lattice isomorphisms between such spaces are unital algebra isomorphisms, which are known \cite{GJ04,AlShakarchiDoust2021} or elementary, and this yields an explicit description of the surjective Hilbert isometries. It turns out that Hilbert's metric distinguishes these spaces. The projective cones of the spaces of H\"{o}lder continuous functions with different exponents, of absolutely continuous functions, and of continuous functions of bounded variation on a compact interval all have the same metric completion, but no two of them are isometric. We also show that the isometries of these projective cones are completely determined by their values on three specifically chosen rays.

The paper is organized as follows. Section~\ref{sec:prelim} contains the preliminaries on order unit spaces, Kakutani's representation theorem, $f$-algebras, and Hilbert's metric. Bounded inverse closed $f$-algebras are studied in Section~\ref{sec:bic}, the main results are proved in Section~\ref{sec:main}, and the applications to function lattices are given in Section~\ref{sec:applications}.

\section{Preliminaries}\label{sec:prelim}

\subsection{Partially ordered vector spaces, order units, and vector lattices}
Let $X$ be a partially ordered vector space and let $X_+$ be its convex cone that determines the order on $X$. The cone $X_+$ is called \emph{Archimedean} if $x\leq 0$ whenever $nx\leq y$ for some $y\in X$ and all $n\in\mathbb N$, and $u\in X_+$ is called a \emph{strong order unit} if for every $x\in X$ there exists $\lambda\geq0$ such that $-\lambda u\leq x\leq\lambda u$. If $X_+$ is an Archimedean cone and $u$ is a strong order unit in $X$, the triple $(X,X_+,u)$ is called an \emph{order unit space}. Every order unit space $(X,X_+,u)$ admits the so-called \emph{order unit norm} $\|\cdot\|_u$ defined by 
\begin{equation*}\label{eq:order-unit-norm}
\|x\|_u:=\inf\{\lambda\geq 0:\; -\lambda u \leq x\leq \lambda u\}.
\end{equation*}
Since $X_+$ is an Archimedean cone, it can be proved that the infimum above is actually attained, and is therefore a minimum. The interior $X_+^\circ$ of the positive cone of the order unit space $(X,X_+,u)$ is always nonempty. In fact, by \cite[Lemma~2.1]{LemRoeWor18}, its interior $X_+^\circ$ coincides with the set of all strong order units of $X$. 

A linear operator $T\colon X\to Y$ between ordered vector spaces is positive if $Tx \in Y_+$ whenever $x \in X_+$. Furthermore, it is readily verified that every positive operator $T\colon (X,X_+,u)\to (Y,Y_+,v)$ between order unit spaces is automatically bounded with $\|T\|=\|Tu\|_v$.

A vector lattice $E$ is a partially ordered vector space such that any two elements of $E$ have a supremum, or equivalently, an infimum. A linear operator $T\colon E\to F$ between vector lattices is a lattice homomorphism if $T(x\vee y)=Tx\vee Ty$, or equivalently, $T(x\wedge y)=Tx\wedge Ty$ for all $x,y \in E$. A bijective lattice homomorphism is called a \emph{lattice isomorphism}; its inverse is again a lattice homomorphism. If $E=F$, a lattice isomorphism is called a \emph{lattice automorphism}. For $e\in E_+$, we denote by $I_e:=\{x\in E:\ -\lambda e\leq x\leq\lambda e \text{ for some }\lambda\geq0\}$ the ideal generated by $e$.

\subsection{Kakutani's representation theorem} 
Let $E$ be an Archimedean vector lattice with a strong order unit $e$. The triple $(E,E_+,e)$ is a special example of an order unit space. Moreover, it can be proved that the order unit norm $\|\cdot\|_e$ is a lattice norm, so that $(E,\|\cdot\|_e)$ is a normed lattice. Its norm completion is an AM-space with a strong order unit $e$. By Kakutani's representation theorem \cite[Theorem~4.29]{AB06}, there exists a compact Hausdorff space $K$ such that $(E,\|\cdot\|_e)$ is lattice isometric to $C(K)$ equipped with the supremum norm. Moreover, the isomorphism can be chosen in such a way that it maps $e$ to the constant one function $\mathbb 1_K$. For future reference, we state this result explicitly. 

\begin{theorem}[Kakutani's representation theorem]\label{thm:Kakutani-representation}
Let $E$ be an Archimedean vector lattice with strong order unit $e$. Then there exist a compact Hausdorff space $K$, unique up to homeomorphism, and a lattice isometry $\Phi\colon (E,\|\cdot\|_e)\to (C(K),\|\cdot\|_\infty)$ such that $\Phi(E)$ is dense in $C(K)$ and $\Phi e=\mathbb 1_K$.
\end{theorem}

Although \cite[Theorem~4.29]{AB06} does not explicitly mention that $e$ is mapped to $\mathbb 1_K$, this is evident from the proof.

Throughout the paper, whenever an Archimedean vector lattice $E$ with strong order unit $e$ is regarded as a norm dense vector sublattice of $C(K)$ for some compact Hausdorff space $K$, this is understood via Kakutani's representation theorem. In particular, $E$ is equipped with the order unit norm $\|\cdot\|_e$, the embedding is a lattice isometry, and $e$ is identified with the constant function $\mathbb 1_K$. We denote the corresponding isometric copy of $E$ in $C(K)$ by $E_K$.

A compact Hausdorff space $K_E$ is called a \emph{Kakutani representation space}, or simply a \emph{Kakutani space}, of $(E,e)$ if there exists a lattice isometric embedding of $(E,\|\cdot\|_e)$ onto a norm dense vector sublattice $E_{K_E}$ of $C(K_E)$ that maps $e$ to $\mathbb 1_{K_E}$. The vector lattice $E_{K_E}$ is called the \emph{Kakutani representation} of $E$. If $K$ is a compact Hausdorff space and $E$ is a norm dense sublattice of $C(K)$ containing the constant one function, then clearly $K$ is the Kakutani space of $E$. In this case, the map $\Phi$ from Theorem~\ref{thm:Kakutani-representation} is the inclusion mapping of $E$ into $C(K)$.

\subsection{$f$-algebras}
An \emph{$f$-algebra} is a vector lattice $E$ equipped with an associative algebra multiplication such that $xy\geq 0$ for all $x,y\geq 0$ and 
\[
xz\wedge y=0 \qquad\text{and}\qquad zx\wedge y=0
\]
for all $x,y\in E$ with $x\wedge y=0$ and all $z\in E_+$. 
If, in addition, there exists an element $e\in E$ satisfying $ex=xe=x$ for every $x\in E$, then $E$ is called a \emph{unital $f$-algebra}, and $e$ is called the \emph{unit} of $E$. 
By the Amemiya--Birkhoff--Pierce theorem \cite[Theorem~2.56]{AB06}, every Archimedean $f$-algebra is commutative. The following lemma will be useful and states that inversion of positive invertible elements in an $f$-algebra with unit is order reversing. 

\begin{lemma}\label{lem:inversion-antimonotone}
    Let $E$ be an $f$-algebra with unit. If $x,y\in E$ are invertible elements satisfying $0\leq x\leq y$, then $0\leq y^{-1}\leq x^{-1}$. 
\end{lemma}
\begin{proof}
We first note that $xy=0$ for all $x,y\in E_+$ with $x\wedge y=0$. Indeed, applying the defining property of an $f$-algebra with $z=y$ gives $xy\wedge y=0$, and applying it once more to the disjoint elements $y$ and $xy$ with $z=x$ gives $xy\wedge xy=0$, that is, $xy=0$. Next, let $x\in E_+$ and $z\in E$. Since $z^+\wedge z^-=0$, we have $xz^-\wedge z^+=0$ and $xz^+\wedge z^-=0$, so that $z^+xz^-=0$ and $z^-xz^+=0$ by the first observation. Hence
\[
zxz=z^+xz^++z^-xz^-\geq0.
\]
In particular, if $x\in E_+$ is invertible, then taking $z=x^{-1}$ gives $x^{-1}=x^{-1}xx^{-1}\geq0$. Now let $x$ and $y$ be invertible with $0\leq x\leq y$. Then $x^{-1},y^{-1}\geq0$, and therefore
\[
x^{-1}-y^{-1}=x^{-1}(y-x)y^{-1}\geq0.\qedhere
\]
\end{proof}

\subsection{Hilbert's metric}\label{subsec:hilbert}

Let $(X,X_+,u)$ be an order unit space with $X\neq\{0\}$. Since $X_+^\circ$ coincides with the set of strong order units of $X$, for all $x,y\in X_+^\circ$ there exist $\alpha,\beta>0$ such that $x \le \alpha y$ and $y\leq \beta x$. Hence
\[
M(x/y):=\inf\{\lambda>0:\; x\leq\lambda y\}
\]
is well-defined. Moreover, $x\leq M(x/y)\,y$, since $X_+$ is Archimedean. Combining $x\leq M(x/y)\,y$ with $y\leq M(y/x)\,x$ gives $x\leq M(x/y)M(y/x)\,x$, and therefore $M(x/y)M(y/x)\geq1$. \emph{Hilbert's metric} is defined by
\[
d_H(x,y):=\log\bigl(M(x/y)\,M(y/x)\bigr), \qquad x,y\in X_+^\circ.
\]
Since $M(\alpha x/\beta y)=\frac{\alpha}{\beta}M(x/y)$ for all $\alpha,\beta>0$, we have $d_H(\alpha x,\beta y)=d_H(x,y)$, so that $d_H(x,y)$ only depends on the rays $[x]:=\{\lambda x:\lambda>0\}$ and $[y]$. Moreover, $d_H$ is symmetric, it satisfies the triangle inequality because $M(x/z)\leq M(x/y)M(y/z)$, and $d_H(x,y)=0$ if and only if $x=\lambda y$ for some $\lambda>0$. Hence $d_H$ is a well-defined metric on the set $P(X_+^\circ):=\{[x]:x\in X_+^\circ\}$ of rays, which we call the \emph{projective cone} of $X$; see also \cite{Nus88}. A \emph{Hilbert isometry} is an isometry with respect to $d_H$.

For a compact Hausdorff space $K$ we have $M(f/g)=\max_K f/g$ for $f,g\in C(K)_+^\circ$. We write $\overline{C(K)}:=C(K)/\mathbb R\one_K$ and equip it with the variation norm
\[
\|\overline h\|_{\mathrm{var}}:=\max_K h-\min_K h.
\] By \cite[Proposition~1.7]{Nus88}, the map
\[
\Log\colon \bigl(P(C(K)_+^\circ),d_H\bigr)\to\bigl(\overline{C(K)},\|\cdot\|_{\mathrm{var}}\bigr),\qquad \Log([f]):=\overline{\log f},
\]
is a surjective isometry; that is, $d_H([f],[g])=\|\overline{\log f-\log g}\|_{\mathrm{var}}$ for all $f,g\in C(K)_+^\circ$. Moreover, by \mbox{\cite[Lemma~4.1]{LemRoeWor18}} the variation norm $\|\cdot\|_{\mathrm{var}}$ is twice the quotient norm of $\|\cdot\|_\infty$ on $\overline{C(K)}$ and the quotient of the Banach space $C(K)$ by the closed subspace $\mathbb R\one_K$ is a Banach space. Hence $(P(C(K)_+^\circ), d_H)$ is a complete metric space. 

The following elementary lemma shows that positive linear bijections with positive inverses, and in particular lattice isomorphisms, induce Hilbert isometries.

\begin{lemma}\label{lem:lattice-isomorphism-order-units-isometry}\label{lem:lattice_isomorphism_Hilbert_isometry}
    Let $T\colon E\to F$ be a positive linear bijection with positive inverse between Archimedean partially ordered vector spaces with strong order units. Then the following statements hold. 
    \begin{enumerate}
        \item $T(E_+^\circ)=F_+^\circ$.
        \item If $u$ and $v$ are strong order units for $E$ and $F$, respectively, and if $Tu=v$, then $T\colon (E,\|\cdot\|_u)\to (F,\|\cdot\|_v)$ is an isometry.
        \item $M(Tx/Ty)=M(x/y)$ for all $x,y\in E_+^\circ$. Consequently, the map $[x]\mapsto[Tx]$ is a well-defined surjective Hilbert isometry from $P(E_+^\circ)$ onto $P(F_+^\circ)$.
    \end{enumerate}  
\end{lemma}

\begin{proof}
    (i). Since $T$ is a positive linear bijection with positive inverse, $x\in E$ is a strong order unit in $E$ if and only if $Tx$ is a strong order unit in $F$. Hence, $T$ maps the set of all strong order units of $E$ onto the set of strong order units of $F$.  

    (ii). Since $T$ is a positive linear bijection with positive inverse, $v=Tu$ is a strong order unit in $F$ and $-\lambda u \leq x\leq \lambda u$ holds if and only if 
    $-\lambda Tu \leq Tx\leq \lambda Tu$ holds. Hence, by the definition of order unit norm, we have 
    $\|Tx\|_{v}=\|x\|_u$. 

    (iii). For $x,y\in E_+^\circ$ and $\lambda>0$ we have $Tx\leq\lambda Ty$ if and only if $x\leq\lambda y$, so $M(Tx/Ty)=M(x/y)$, and therefore $d_H([Tx],[Ty])=d_H([x],[y])$. Since $T$ is linear, $[x]\mapsto[Tx]$ is well defined, and it is surjective by~(i).
\end{proof}

In the case of $f$-algebras, inversion of positive invertible elements reverses the order in the following multiplicative sense. Consequently, inversion induces an isometry for Hilbert's metric.

\begin{lemma}\label{lem:inversion-order-reversing}
Let $E$ be an Archimedean $f$-algebra with unit. If the unit of $E$ is a strong order unit and $x,y\in E_+^\circ$ are invertible with inverses in $E_+^\circ$, then
\begin{equation*}
M(x^{-1}/y^{-1})=M(y/x) \qquad\text{and}\qquad d_H\bigl([x^{-1}],[y^{-1}]\bigr)=d_H\bigl([x],[y]\bigr).
\end{equation*}
\end{lemma}

\begin{proof}
It follows from Lemma~\ref{lem:inversion-antimonotone} that if $x,y\in E_+$ are invertible, then for every $\lambda>0$ we have $x^{-1}\leq\lambda y^{-1}$ if and only if $y\leq\lambda x$. Hence \[\{\lambda>0: x^{-1}\leq\lambda y^{-1}\}=\{\lambda>0: y\leq\lambda x\},\] so $M(x^{-1}/y^{-1})=M(y/x)$. Interchanging $x$ and $y$ gives $M(y^{-1}/x^{-1})=M(x/y)$, and the identity for $d_H$ follows from the definition.
\end{proof}

In particular, if the unit of $E$ is a strong order unit and every $x\in E_+^\circ$ is invertible with $x^{-1}\in E_+^\circ$, then the \emph{inversion} $[x]\mapsto[x^{-1}]$ is a well-defined map on $P(E_+^\circ)$, since $(\lambda x)^{-1}=\lambda^{-1}x^{-1}$ for all $\lambda>0$. By Lemma~\ref{lem:inversion-order-reversing}, it is a surjective Hilbert isometry of $P(E_+^\circ)$, which is its own inverse.

The results in this paper rely on the following description of the surjective Hilbert isometries between projective cones of spaces of continuous functions. In Section~\ref{sec:main} we extend it to arbitrary Archimedean vector lattices with strong order units.

\begin{theorem}[{\cite[Theorem~4.7]{LemRoeWor18}}]\label{thm:Banach-Stone-projective}
Let $K$ and $L$ be compact Hausdorff spaces. A map $\psi\colon P(C(K)_+^\circ)\to P(C(L)_+^\circ)$ is a surjective Hilbert isometry if and only if there exist $g\in C(L)_+^\circ$, $\varepsilon\in\{\pm1\}$, and a homeomorphism $\vartheta\colon L\to K$ such that
\[
\psi([f])=[g\,(f\circ\vartheta)^\varepsilon]\qquad\text{for all } [f]\in P(C(K)_+^\circ).
\]
\end{theorem}
Note that the reverse implication in Theorem~\ref{thm:Banach-Stone-projective} follows from Lemmas~\ref{lem:lattice-isomorphism-order-units-isometry} and~\ref{lem:inversion-order-reversing}. Indeed, $f\mapsto g\,(f\circ\vartheta)$ is a lattice isomorphism of $C(K)$ onto $C(L)$, so $[f]\mapsto[g\,(f\circ\vartheta)]$ is a surjective Hilbert isometry by Lemma~\ref{lem:lattice-isomorphism-order-units-isometry}, and the inversion $[f]\mapsto[f^{-1}]$ is a surjective Hilbert isometry of $P(C(L)_+^\circ)$ by Lemma~\ref{lem:inversion-order-reversing}.

\section{Bounded inverse closed $f$-algebras}\label{sec:bic}

By Theorem~\ref{thm:Banach-Stone-projective}, the surjective Hilbert isometries of $P(C(K)_+^\circ)$ may involve the inversion $[f]\mapsto[f^{-1}]$, which uses the multiplicative structure of $C(K)$. For a general Archimedean vector lattice with a strong order unit, no algebraic structure is available a priori. In this section we introduce the class of bounded inverse closed $f$-algebras and show that being part of this class can be read off from any Kakutani representation. Furthermore, we give criteria that will be used in Section~\ref{sec:main} to help decide when this inversion is actually available.

\begin{definition}\label{def:bic}
Let $E$ be a unital $f$-algebra with unit $e$.
\begin{enumerate}
\item $E$ is said to be \emph{bounded inverse closed} if every $x\in E$ with $x\geq e$ is invertible in $E$.
\item A subalgebra $F$ of $E$ containing $e$ is said to be \emph{bounded inverse closed} (in $E$) if, whenever $x\in F$ satisfies $x\geq e$ and is invertible in $E$, its inverse $x^{-1}$ belongs to $F$.
\end{enumerate}
\end{definition}

Note that we do not assume $E$ to be Archimedean in Definition~\ref{def:bic}; whenever this assumption is needed, it will be mentioned explicitly. Property~(i) in Definition~\ref{def:bic} is known in the theory of $f$-rings as the \emph{bounded inversion property}; see, for instance, \cite{HJ61,HIJ61}.

For Archimedean vector lattices with a strong order unit, the following theorem shows that both properties in Definition~\ref{def:bic} are determined by the Kakutani representation.

\begin{theorem}\label{thm:f-algebra-structure-kakutani-representation}
Let $E$ be an Archimedean vector lattice with strong order unit $u$, let
$K_E$ be a Kakutani space of $E$, and let
$\Phi_E\colon E\to C(K_E)$
be the corresponding lattice isometry with dense range satisfying
$\Phi_E(u)=\mathbb 1_{K_E}$. Then the following statements hold.
\begin{enumerate}
\item
The vector lattice $E$ admits a multiplication that turns it into a
unital $f$-algebra with unit $u$ if and only if $E_{K_E}$ is a
subalgebra of $C(K_E)$.

\item
The vector lattice $E$ admits a multiplication that turns it into a
bounded inverse closed $f$-algebra with unit $u$ if and only if
$E_{K_E}$ is a bounded inverse closed subalgebra of $C(K_E)$.
\end{enumerate}
In either case, the multiplication is unique and given by
\[
x\cdot y
=
\Phi_E^{-1}\bigl(\Phi_E(x)\Phi_E(y)\bigr)
\]
for all $x,y\in E$, and it is independent of the chosen Kakutani
representation.
\end{theorem}

\begin{proof}
(i). Suppose first that $E_{K_E}$ is a subalgebra of $C(K_E)$. For
$x,y\in E$, define
\begin{equation*}
x\cdot y
:=
\Phi_E^{-1}\big(\Phi_E(x)\Phi_E(y)\big).
\end{equation*}
Since $E_{K_E}$ is closed under pointwise multiplication, this operation is
well defined. The associativity, commutativity, and distributivity of this
multiplication follow immediately from the corresponding properties of
pointwise multiplication in $C(K_E)$. Moreover, since $\Phi_E(u)=\mathbb 1_{K_E}$ and $\Phi_E(x\cdot y)=\Phi_E(x)\Phi_E(y)$ for all $x,y\in E$, $u$ is the multiplicative unit of $E$ for the product defined above.

Now we prove that $E$ equipped with this product is an $f$-algebra. To that end, pick $x,y\in E_+$ and note that $\Phi_E(x)$ and $\Phi_E(y)$ are nonnegative functions. This yields 
\[
\Phi_E(x\cdot y)=\Phi_E(x)\Phi_E(y)\geq0,
\]
and since $\Phi_E\colon E\to E_{K_E}$ is a lattice isomorphism, we have $x\cdot y\geq 0$ in $E$. Furthermore, if $x,y\in E$ satisfy
$x\wedge y=0$, then $\Phi_E(x)\wedge\Phi_E(y)=0$. Hence, for any $z\in E_+$ the pointwise multiplication implies
\[
\Phi_E(x)\Phi_E(z)\wedge\Phi_E(y)=0,
\]
and so, $(x\cdot z)\wedge y = (z\cdot x)\wedge y = 0$. This proves that the multiplication defined above turns $E$ into a unital $f$-algebra with unit $u$.

Conversely, suppose that $E$ is equipped with a multiplication, denoted by juxtaposition, that turns it into a unital $f$-algebra with unit $u$. Then $E$ and $C(K_E)$ are Archimedean $f$-algebras with unit elements $u$ and $\mathbb 1_{K_E}$, respectively, and $\Phi_E$ is a lattice homomorphism, in particular a positive linear map, with $\Phi_E(u)=\mathbb 1_{K_E}$. By \cite[Corollary~5.5]{HdP84}, $\Phi_E$ is an algebra homomorphism, so that
\[
\Phi_E(x)\Phi_E(y)=\Phi_E(xy)\in E_{K_E}\qquad\text{for all }x,y\in E.
\]
Hence $E_{K_E}$ is a subalgebra of $C(K_E)$, and $xy=\Phi_E^{-1}\bigl(\Phi_E(x)\Phi_E(y)\bigr)$ for all $x,y\in E$. In particular, the multiplication is unique. Since this argument applies to every Kakutani representation of $E$, the multiplication is also independent of the chosen Kakutani representation.

(ii). Suppose now that $E_{K_E}$ is a bounded inverse closed subalgebra of
$C(K_E)$, and equip $E$ with the multiplication induced by
$\Phi_E$. Let $x\in E$ satisfy $x\geq u$. Then
$\Phi_E(x)\geq\mathbb 1_{K_E}$. Since $E_{K_E}$ is bounded inverse closed, we have $\Phi_E(x)^{-1}\in E_{K_E}$. The element $y:=\Phi_E^{-1}\big(\Phi_E(x)^{-1}\big)$ satisfies 
\[
\Phi_E(x\cdot y)=\Phi_E(x)\Phi_E(x)^{-1}=\mathbb 1_{K_E}=\Phi_E(u),
\]
and since $\Phi_E$ is injective, we have $x\cdot y = y \cdot x = u$ since $E$ is commutative. It follows that every
$x\geq u$ is invertible in $E$. Hence, $E$ is bounded inverse closed. 

Conversely, suppose that $E$ is a bounded inverse closed $f$-algebra with unit $u$. By (i) and its proof, $E_{K_E}$ is a subalgebra of $C(K_E)$ and $\Phi_E\colon E\to E_{K_E}$ is an algebra isomorphism. Let $f\in E_{K_E}$ satisfy
$f\geq\mathbb 1_{K_E}$ and write $f=\Phi_E(x)$ for some $x\in E$.
Since $\Phi_E$ is an order isomorphism, we have $x\geq u$. Hence $x$ has an
inverse $x^{-1}\in E$. It follows that
\[
f^{-1}=\Phi_E(x)^{-1}=\Phi_E(x^{-1})\in E_{K_E},
\]
so that $E_{K_E}$ is bounded inverse closed.
\end{proof} 

A lattice isomorphism that maps the order unit to the order unit can be used to transport a multiplication from one vector lattice to the other. Consequently, these properties are invariant under such isomorphisms.

\begin{corollary}\label{cor:invariance-algebraic-kakutani-properties}
Let $E$ and $F$ be Archimedean vector lattices with strong order units
$u\in E_+$ and $v\in F_+$, respectively, and let
$T\colon E\to F$ be a lattice isomorphism satisfying $Tu=v$.
Let $K_E$ and $K_F$ be Kakutani spaces of $E$ and $F$, respectively.
Then the following statements hold.
\begin{enumerate}
\item
The sublattice $E_{K_E}$ is a subalgebra of $C(K_E)$ if and only if
$F_{K_F}$ is a subalgebra of $C(K_F)$.

\item
The sublattice $E_{K_E}$ is a bounded inverse closed subalgebra of
$C(K_E)$ if and only if $F_{K_F}$ is a bounded inverse closed subalgebra
of $C(K_F)$.

\end{enumerate}
\end{corollary}
 
\begin{proof}
(i). Suppose that $E_{K_E}$ is a subalgebra of $C(K_E)$. By
Theorem~\ref{thm:f-algebra-structure-kakutani-representation}, $E$ admits a multiplication that turns it into a unital $f$-algebra with unit $u$.
Define a multiplication on $F$ by 
\[
x\star y
:=
T\big(T^{-1}x\,T^{-1}y\big)
\] 
for all $x,y\in F$. Since $T$ is a lattice isomorphism and $Tu=v$, this multiplication turns
$F$ into a unital $f$-algebra with unit $v$. Another application of
Theorem~\ref{thm:f-algebra-structure-kakutani-representation} shows that
$F_{K_F}$ is a subalgebra of $C(K_F)$. The converse follows by applying
the same argument to $T^{-1}$.

(ii). If $E_{K_E}$ is bounded inverse closed, then by Theorem~\ref{thm:f-algebra-structure-kakutani-representation} we have that $E$ is a bounded inverse closed $f$-algebra. The multiplication
$\star$ defined in (i) has the same property on $F$, because $x\geq v$ implies
$T^{-1}x\geq u$, and the inverse of $x$ is the image under $T$ of the
inverse of $T^{-1}x$. Hence $F_{K_F}$ is bounded inverse closed by
Theorem~\ref{thm:f-algebra-structure-kakutani-representation}. The reverse
implication follows analogously.  
\end{proof}

Corollary~\ref{cor:invariance-algebraic-kakutani-properties} shows that whether $E_{K_E}$ is a subalgebra, or a bounded inverse closed subalgebra, of $C(K_E)$ does not depend on the choice of the Kakutani representation of $E$. We shall therefore say that an Archimedean vector lattice with a strong order unit is a unital $f$-algebra (respectively, a bounded inverse closed $f$-algebra) if one, and hence every, Kakutani representation has the corresponding property. In this case, the multiplication on the vector lattice is the one induced by the pointwise multiplication on any Kakutani representation.

In concrete situations, it is often easier to check that a subspace is closed under inversion than that it is a subalgebra. The following proposition shows that, for subspaces containing the unit, being closed under taking inverses in this sense already forces the subspace to be closed under multiplication.

\begin{proposition}\label{prop:bounded-inverse-closed-subspace-subalgebra}
Let $E$ be a commutative $f$-algebra with unit $e$, and let $V$ be a vector subspace of $I_e$ containing $e$. If every $x\in V$ satisfying $x\geq e$ is invertible and satisfies $x^{-1}\in V$, then $V$ is a bounded inverse closed subalgebra of $E$.
\end{proposition}

\begin{proof}
We first observe that if $x\in V$ and there exists $\alpha>0$ such that $\alpha x\geq e$, then $x^{-1}\in V$. We show that $x^2\in V$ for every $x\in V$. To this end, choose $x\in V$. Since $x\in I_e$, there exists $\lambda>0$ such that $-\lambda e\leq x\leq \lambda e$. Increasing $\lambda$ if necessary, we may assume that $y:=x+\lambda e\geq e$ and $y+e\geq e$. By assumption, $y$ and $y+e$ are invertible, and their inverses $y^{-1}$ and $(y+e)^{-1}$ are elements of $V$. Therefore,
\[
z := y^{-1}-(y+e)^{-1}=(y(y+e))^{-1}
\]
belongs to $V$. Since $y\in I_e$, there exists $\mu>0$ such that $e\leq y\leq \mu e$, which yields $e\leq y+e\leq (\mu+1)e$. Hence
\[
e\leq y(y+e)\leq \mu(\mu+1)e,
\]
so that $\mu(\mu+1)(y(y+e))^{-1}\geq e$ by Lemma~\ref{lem:inversion-antimonotone}. It follows that $z^{-1}\in V$, and so
\[
y^2+y=y(y+e)=z^{-1}\in V.
\]
Since $y\in V$, we have $y^2\in V$. As $x=y-\lambda e$, we obtain
\[
x^2=y^2-2\lambda y+\lambda^2e\in V.
\]
Since $E$ is commutative, we have $xy=\frac14\big((x+y)^2-(x-y)^2\big)\in V$ for all $x,y\in V$. Hence $V$ is a subalgebra, which is bounded inverse closed.
\end{proof}

The following corollary is a special case of Proposition~\ref{prop:bounded-inverse-closed-subspace-subalgebra}.

\begin{corollary}\label{cor:bic-subspace-CK}
Let $K$ be a compact Hausdorff space and let $V$ be a vector subspace
of $C(K)$ containing $\mathbb 1$. Suppose that $f^{-1}\in V$ for every
$f\in V$ with $f\geq \mathbb 1$. Then $V$ is a subalgebra of $C(K)$.
\end{corollary}

For subspaces of $C(K)$, the following corollary reformulates this in terms of the interior of the positive cone. This is the form used in Section~\ref{sec:main}.

\begin{corollary}\label{lem:inverse-isometry-well-defined}
Let $K$ be a compact Hausdorff space and let $V$ be a vector subspace of $C(K)$ containing $\mathbb 1$. Then $f^{-1}\in V$ for every $f\in V\cap C(K)_+^\circ$ if and only if $V$ is a bounded inverse closed subalgebra of $C(K)$. In this case, $f^{-1}\in V\cap C(K)_+^\circ$ for every $f\in V\cap C(K)_+^\circ$.
\end{corollary}

\begin{proof}
If $f^{-1}\in V$ for every $f\in V\cap C(K)_+^\circ$, then in particular $f^{-1}\in V$ for every $f\in V$ with $f\geq\mathbb 1$. Hence $V$ is a subalgebra of $C(K)$ by Corollary~\ref{cor:bic-subspace-CK}, and it is bounded inverse closed by definition. Conversely, suppose that $V$ is a bounded inverse closed subalgebra of $C(K)$ and let $f\in V\cap C(K)_+^\circ$. Since $K$ is compact, there are $\delta,\lambda>0$ such that $\delta\mathbb 1\leq f\leq\lambda\mathbb 1$. Since $\delta^{-1}f\geq\mathbb 1$, we have $\delta f^{-1}=(\delta^{-1}f)^{-1}\in V$, and hence $f^{-1}\in V$. Moreover, $f^{-1}\geq\lambda^{-1}\mathbb 1$, so $f^{-1}\in C(K)_+^\circ$.
\end{proof}

Note that if $V$ is a norm dense sublattice of $C(K)$ with $\mathbb 1\in V$, as is the case for a Kakutani representation, then $V\cap C(K)_+^\circ=V_+^\circ$. In this case, Corollary~\ref{lem:inverse-isometry-well-defined} shows that inversion maps $V_+^\circ$ into itself precisely when $V$ is a bounded inverse closed subalgebra of $C(K)$.

\begin{example}\label{ex:prop-sharp}
The following three examples show that the assumptions in Proposition~\ref{prop:bounded-inverse-closed-subspace-subalgebra} cannot be weakened.
\begin{enumerate}
\item It is essential that every $x\in V$ with $x\geq e$ is invertible. Let $E:=\mathbb R[t]$, ordered by setting a nonzero polynomial $p(t)=a_0+a_1t+\cdots+a_nt^n$ to be positive if $a_k>0$, where $k:=\min\{i:a_i\neq0\}$. Then $E$ is a totally ordered, non-Archimedean commutative $f$-algebra with strong order unit $\mathbb 1$. In $V:=\Span\{\mathbb 1,t\}$, the invertible elements are the nonzero constants, whose inverses belong to $V$. However, $2\cdot\mathbb 1+t\geq\mathbb 1$ is not invertible, and $V$ is not a subalgebra, since $t^2\notin V$.
\item The assumption $V\subseteq I_e$ cannot be omitted. In $E:=C(\mathbb R)$, the subspace $V:=\Span\{\mathbb 1,\mathrm{id}\}$ is not contained in $I_{\mathbb 1}=C_b(\mathbb R)$. If $f=\alpha\mathbb 1+\beta\,\mathrm{id}\geq\mathbb 1$, then $\beta=0$ and $\alpha\geq1$, so that $f^{-1}\in V$. However, $\mathrm{id}^2\notin V$.
\item The assumption $e\in V$ cannot be omitted. The subspace $V:=\Span\{\mathrm{id}\}$ of $C[0,1]$ contains no $f\geq\mathbb 1$, so that the hypothesis of Proposition~\ref{prop:bounded-inverse-closed-subspace-subalgebra} holds vacuously, but $\mathrm{id}^2\notin V$.
\end{enumerate}
\end{example}

\section{Hilbert isometries of Archimedean vector lattices with strong order units}\label{sec:main}

In this section we describe the surjective Hilbert isometries between projective cones of Archimedean vector lattices with strong order units. In Subsection~\ref{subsec:extension} we show that Hilbert isometries between the projective cones of such sublattices extend uniquely to the projective cones of the corresponding spaces of continuous functions, so that Theorem~\ref{thm:Banach-Stone-projective} applies to the extension. The main difficulty is to show that the inversion can only occur when the vector lattices are bounded inverse closed $f$-algebras, which is done in Subsection~\ref{subsec:eps-minus-one}. The main theorem is proved in Subsection~\ref{subsec:main-theorem}, and the isometry group is determined in Subsection~\ref{subsec:isometry-group}.

\subsection{Extending isometries to $P(C(K)_+^\circ)$}\label{subsec:extension}

We first show that the projective cone of a norm dense sublattice is dense in the projective cone of the whole space.

\begin{lemma}\label{lem:dense_sublattice_dense_projectice_cone}
Let $K$ be a compact Hausdorff space and let $E$ be a norm dense sublattice of $C(K)$ that contains the constant one function. Then $P(E_+^\circ)$ is dense in $P(C(K)_+^\circ)$ for Hilbert's metric.    
\end{lemma}

\begin{proof}
    We claim that $E_+^\circ$ is norm dense in $C(K)_+^\circ$. 
    Let $f\in C(K)_+^\circ$ be arbitrary. Then there exists $\delta >0$ such that $f\geq 2\delta \one$. Pick $0<\varepsilon<\delta$. Since $E$ is dense in $C(K)$, there exists $g\in E$ such that $\|f-g\|_\infty<\varepsilon$. From $|f-g|<\varepsilon \one$, we obtain 
\[
g\geq f-\varepsilon \one \geq f-\delta\mathbb 1\geq \delta \mathbb 1
\]
 which shows $g\in E_+^\circ$. 

    To prove that $P(E_+^\circ)$ is dense in $P(C(K)_+^\circ)$, take any $[f]\in P(C(K)_+^\circ)$. Then for each $n\in \mathbb N$ there exists $g_n\in E_+^\circ$ such that $\|f-g_n\|_\infty < \tfrac{1}{n}$. We claim that 
    $\overline{g_n}\to\overline{f}$ with respect to the variation norm on $\overline{C(K)}$. To this end, since $f\in C(K)_+^\circ$, there exists $\varepsilon>0$ such that $f\geq 2\varepsilon \one$. Moreover, from the convergence $g_n\to f$ in $C(K)$ we conclude that there exists $N\in \mathbb N$ such that $g_n\geq \varepsilon\one$ for all $n\geq N$. For $t \in K$, note that the mean value theorem yields 
\[
|\log f(t) - \log g_n(t)| = c^{-1}|f(t) - g_n(t)| \le c^{-1}\|f - g_n\|_\infty
\]
 for some $c$ between $f(t)$ and $g_n(t)$. Since $f(t), g_n(t) \ge \varepsilon$, we obtain  
    
\[
\|\log f-\log g_n\|_\infty \leq \frac{1}{\varepsilon}\|f-g_n\|_\infty
\]
 for all $n\geq N$. By \cite[Proposition~1.7]{Nus88} the mapping $\Log\colon P(C(K)_+^\circ) \to \overline{C(K)}$ is an isometry, and we have that $\|\overline h\|_{\mathrm{var}} \leq 2\|h\|_\infty$ for every $h\in C(K)$, so it follows that
    \begin{align*}
    d_H([f],[g_n])&=\|\Log([f])-\Log([g_n])\|_{\mathrm{var}}=\|\overline{\log f-\log g_n}\|_{\mathrm{var}}\\&\leq 2\|\log f-\log g_n\|_\infty
    \leq \frac{2}{\varepsilon}\|f-g_n\|_\infty\to 0.\qedhere
    \end{align*}
\end{proof}

Since the projective cone of $C(K)$ is complete for Hilbert's metric, Lemma~\ref{lem:dense_sublattice_dense_projectice_cone} allows us to extend maps between projective cones of dense sublattices.

\begin{lemma}\label{L: unique extension of h}
    Let $K$, $L$ be compact Hausdorff spaces and let $E\subseteq C(K)$ and $F\subseteq C(L)$ be norm dense vector sublattices containing the constant one function. Further, suppose that we have a map
\[
\tau \colon P(E_+^\circ)\to P(F_+^\circ).
\]
For Hilbert's metric the following statements hold.
\begin{enumerate}
\item 
If $\tau$ is uniformly continuous, then $\tau$ admits a unique uniformly continuous extension
\[
\overline \tau\colon
P(C(K)_+^\circ)\to P(C(L)_+^\circ).
\]
\item 
If $\tau$ is Lipschitz with constant $\alpha\geq0$, then $\overline{\tau}$ is Lipschitz with constant $\alpha$.
\item
If $\tau$ is an isometric embedding, then $\overline \tau$ is an isometric embedding.
\item
If $\tau$ is a surjective isometry, then $\overline \tau$ is a surjective isometry.
\end{enumerate}
\end{lemma}

\begin{proof}
Since $\one_K$ is a strong order unit of $E$, we have $E_+^\circ = E\cap C(L)_+^\circ$, and $M(f/g)$ is the same whether it is computed in $E$ or in $C(K)$. Hence $(P(E_+^\circ),d_H)$ is a metric subspace of $(P(C(K)_+^\circ),d_H)$, and similarly for $F$ and $C(L)$. As mentioned in the preliminaries, the metric space $(P(C(K)_+^\circ),d_H)$ is complete. Furthermore, $P(E_+^\circ)$ is dense in $P(C(K)_+^\circ)$ and $P(F_+^\circ)$ is dense in $P(C(L)_+^\circ)$ by Lemma~\ref{lem:dense_sublattice_dense_projectice_cone}.

(i). The map $\tau$ is uniformly continuous, its domain is dense in $P(C(K)_+^\circ)$, and its codomain is contained in the complete metric space $P(C(L)_+^\circ)$. It therefore follows from \mbox{\cite[Theorem~4.3.17]{Eng89}} that $\tau$ extends to a uniformly continuous map $\overline\tau$ from $P(C(K)_+^\circ)$ into $P(C(L)_+^\circ)$. Any two continuous extensions of $\tau$ agree on the dense set $P(E_+^\circ)$, and hence everywhere by \mbox{\cite[Theorem~2.1.9]{Eng89}}. Thus $\overline\tau$ is the unique continuous extension of $\tau$.

(ii) and (iii). Let $[f],[g]\in P(C(K)_+^\circ)$, and choose sequences $([f_n])$ and $([g_n])$ in $P(E_+^\circ)$ converging to $[f]$ and $[g]$, respectively. Since $\overline\tau$ and $d_H$ are continuous, we have
\begin{align*}
d_H\bigl(\overline\tau([f]),\overline\tau([g])\bigr)&=\lim_{n\to\infty} d_H\bigl(\tau([f_n]),\tau([g_n])\bigr),\\
d_H([f],[g])&=\lim_{n\to\infty} d_H([f_n],[g_n]).
\end{align*}
Applying the Lipschitz estimate for $\tau$ to the pairs $[f_n],[g_n]$ and letting $n\to\infty$ yields (ii). If $\tau$ is an isometric embedding, the same argument with an equality instead of an inequality yields (iii).

(iv). Suppose that $\tau$ is a surjective isometry. The sets $P(E_+^\circ)$ and $P(F_+^\circ)$ are dense in the complete metric spaces $P(C(K)_+^\circ)$ and $P(C(L)_+^\circ)$, respectively. Therefore, by \mbox{\cite[Corollary~4.3.18]{Eng89}}, $\tau$ extends to a surjective isometry between these two spaces, and by the uniqueness in (i) this extension is $\overline\tau$.
\end{proof}

\subsection{The case $\varepsilon=-1$}\label{subsec:eps-minus-one}
Let $E$ and $F$ be Archimedean vector lattices with strong order units, and let $\psi  \colon P(E_+^\circ) \to P(F_+^\circ)$ be a surjective Hilbert isometry. By Kakutani's representation theorem, we may regard $E$ and $F$ as norm dense sublattices of $C(K)$ and $C(L)$, respectively, both containing the constant one function. Now $\psi$ extends uniquely to a surjective isometry 
    \[
    \overline \psi  \colon P(C(K)_+^\circ) \to P(C(L)_+^\circ)
    \]
    by Lemma~\ref{L: unique extension of h}. By Theorem~\ref{thm:Banach-Stone-projective} there exist $g \in C(L)_+^\circ$, $\varepsilon \in \{\pm 1\}$, and a homeomorphism $\vartheta \colon L \to K$, such that 
    \[
    \overline{\psi}([f]) = [g (f \circ \vartheta)^\varepsilon] \quad \mbox{for all $[f] \in P(C(K)_+^\circ)$.}
    \]
Since neither $E$ nor $F$ is assumed to be equipped with the structure of multiplication, the case $\varepsilon = -1$ requires some further investigation. To that end, we define the subsets 
\begin{align}\label{E:F_{-1}}
F_{-1} := \{f\circ \vartheta \colon f \in E_+^\circ\}  \quad\mathrm{and}\quad F_0 := \{ f\circ \vartheta \colon f \in E\} = F_{-1} - F_{-1}
\end{align}
of $C(L)$.
We first show that $F_{0}$ is a bounded inverse closed subalgebra of $C(L)$. 

\begin{proposition}\label{P:eps = -1}
    Let $E$ and $F$ be Archimedean vector lattices with strong order units, regarded as norm dense sublattices of $C(K)$ and $C(L)$ containing the constant one function. Suppose that $\psi\colon P(E_+^\circ)\to P(F_+^\circ)$ is a surjective Hilbert isometry whose unique extension $\overline\psi\colon P(C(K)_+^\circ)\to P(C(L)_+^\circ)$ is given by
\[
\overline\psi([f])=[g(f\circ\vartheta)^{-1}]\qquad\text{for all }[f]\in P(C(K)_+^\circ)
\]
for some homeomorphism $\vartheta\colon L\to K$ and some $g\in F_+^\circ$, and let $F_{-1}$ and $F_0$ be as in \eqref{E:F_{-1}}. Then $F_{-1}$ is closed under inversion, and $F_{0}$ is a bounded inverse closed subalgebra of $C(L)$.  
\end{proposition}

\begin{proof}
    Let $k \in F_{-1}$ and define $l := g  k^{-1} + g$. Writing $k=f\circ\vartheta$ with $f\in E_+^\circ$, we have $\psi([f])=[gk^{-1}]\in P(F_+^\circ)$, so $gk^{-1}\in F_+^\circ$. Since also $g\in F_+^\circ$, it follows that $l\in F_+^\circ$. As $\psi$ is surjective, there is $m\in E_+^\circ$ with $\psi([m])=[l]$; replacing $m$ by a suitable positive multiple, we may assume that $g(m\circ\vartheta)^{-1}=l$. Hence $n := m \circ \vartheta \in F_{-1}$ satisfies $g n^{-1} = gk^{-1} + g$. Since $g$ is invertible in $C(L)$, we find that $n = k(k + \mathbb{1}_L)^{-1}$. Hence,
    \[
    (\mathbb 1_K-m)\circ \vartheta =\mathbb 1_L-n=\mathbb 1_L-k (k + \mathbb{1}_L)^{-1}=(k+\mathbb 1_L)^{-1} \in C(L)_+^\circ,
    \]
    which shows that $\mathbb 1_K-m$ is positive and invertible in $C(K)$. Since $\mathbb 1_K-m$ belongs to $E_+^\circ$, we have $(k+\mathbb 1_L)^{-1}=(\mathbb{1}_K - m)\circ \vartheta \in F_{-1}$. 

    With $k=f\circ\vartheta$ as above, there exists $\lambda>0$ such that $\lambda f-\mathbb 1_K\in E_+^\circ$. Since  
    \[
    \lambda k-\mathbb 1_L=\lambda (f\circ \vartheta)-\mathbb 1_K\circ \vartheta=(\lambda f-\mathbb 1_K)\circ \vartheta,
    \] 
    we have $\lambda k-\mathbb 1_L \in F_{-1}$. By the above, we have 
    \[
    \lambda^{-1}k^{-1}=(\lambda k)^{-1}=\big((\lambda k-\mathbb 1_L) + \mathbb 1_L\big)^{-1}\in F_{-1},
    \]
    and therefore, $k^{-1}\in F_{-1}$. Now let $x\in F_0$ satisfy $x\geq\mathbb 1_L$. Then $x=f\circ\vartheta$ for some $f\in E$ with $f\geq\mathbb 1_K$, so that $x\in F_{-1}$, and hence $x^{-1}\in F_{-1}\subseteq F_0$. Since $\mathbb 1_L\in F_0$, it follows from Proposition~\ref{prop:bounded-inverse-closed-subspace-subalgebra} that $F_0$ is a bounded inverse closed subalgebra of $C(L)$.
\end{proof}

\subsection{The main theorem}\label{subsec:main-theorem}
Combining the extension from Subsection~\ref{subsec:extension} with Proposition~\ref{P:eps = -1}, we obtain the following description of the surjective Hilbert isometries.

\begin{theorem}\label{T:d_H-isom between E and F}
    Let $E$ and $F$ be Archimedean vector lattices with strong order units. A map $\psi  \colon P(E_+^\circ) \to P(F_+^\circ)$ is a surjective Hilbert isometry if and only if there exist an $\varepsilon\in \{\pm 1\}$ and a lattice isomorphism $T\colon E\to F$ such that 
    \[
    \psi([f]) = [(Tf)^\varepsilon]
    \] for all $[f] \in P(E_+^\circ)$.
\end{theorem}

\begin{proof}
   \begingroup\relpenalty=10000 \binoppenalty=10000
($\Rightarrow$). Let $\psi$ be a surjective Hilbert isometry. As mentioned in \mbox{Subsection~\ref{subsec:eps-minus-one}}, we regard $E$ and $F$ as norm dense sublattices of $C(K)$ and $C(L)$, respectively, containing the constant one functions, and we obtain $g\in C(L)_+^\circ$, $\varepsilon\in\{\pm1\}$, and a homeomorphism $\vartheta\colon L\to K$ such that
\[
\psi([f])=[g\,(f\circ\vartheta)^\varepsilon]\qquad\text{for all } [f]\in P(E_+^\circ).
\]
Since $\mathbb 1_K\in E_+^\circ$, we have $[g]=\psi([\mathbb 1_K])\in P(F_+^\circ)$, so $g$ is a positive multiple of an element of $F_+^\circ$. Since $F_+^\circ$ is closed under multiplication by positive scalars, it follows that $g\in F_+^\circ$.

Suppose first that $\varepsilon=1$, and let $T\colon C(K)\to C(L)$ be the lattice isomorphism $Tf:=g\,(f\circ\vartheta)$. For $f\in E_+^\circ$ we have $[Tf]=\psi([f])\in P(F_+^\circ)$, so $Tf\in F_+^\circ$. Since $\psi$ is surjective, every element of $F_+^\circ$ is of the form $\lambda Tf=T(\lambda f)$ for some $\lambda>0$ and $f\in E_+^\circ$. Hence $T(E_+^\circ)=F_+^\circ$, and since $E_+^\circ$ and $F_+^\circ$ span $E$ and $F$, respectively, $T$ restricts to a lattice isomorphism from $E$ onto $F$.

Suppose now that $\varepsilon=-1$, and consider the sets $F_0$ and $F_{-1}$ defined in~\eqref{E:F_{-1}}. By Proposition~\ref{P:eps = -1}, $F_{-1}$ is closed under inversion, and it is clearly closed under multiplication by positive scalars. We claim that $F_+^\circ=gF_{-1}$. Let $k\in F_{-1}$, and write $k^{-1}=f\circ\vartheta$ with $f\in E_+^\circ$. Then $[gk]=[g(f\circ\vartheta)^{-1}]=\psi([f])\in P(F_+^\circ)$, so $gk\in F_+^\circ$. This shows that $gF_{-1}\subseteq F_+^\circ$. Conversely, let $x\in F_+^\circ$. Since $\psi$ is surjective, there exist $f\in E_+^\circ$ and $\lambda>0$ such that $x=\lambda g(f\circ\vartheta)^{-1}$. As $\lambda(f\circ\vartheta)^{-1}\in F_{-1}$, we obtain $x\in gF_{-1}$, which proves the claim.

We proceed to show that $F_0 = F$. Since $\mathbb 1_L\in F_+^\circ=gF_{-1}$, there is $k\in F_{-1}$ with $gk=\mathbb 1_L$, so $g^{-1}=k\in F_{-1}$. As $F_{-1}$ is closed under inversion, also $g\in F_{-1}$. Moreover, $F_{-1}=F_0\cap C(L)_+^\circ$, and $F_0$ is a subalgebra of $C(L)$ by Proposition~\ref{P:eps = -1}. As $g,g^{-1}\in F_{-1}$, it follows that multiplication by $g$ maps $F_{-1}$ bijectively onto itself. Therefore
\[
F_{-1}=gF_{-1}=F_+^\circ,
\]
and since $F_{-1}$ and $F_+^\circ$ span $F_0$ and $F$, respectively, we obtain $F_0=F$.

Finally, let $T\colon C(K)\to C(L)$ be the lattice isomorphism $Tf:=g^{-1}(f\circ\vartheta)$. The composition operator $f\mapsto f\circ\vartheta$ maps $E_+^\circ$ onto $F_{-1}=F_+^\circ$, and $g^{-1}F_{-1}=F_{-1}$, so $T(E_+^\circ)=F_+^\circ$. As before, $T$ restricts to a lattice isomorphism from $E$ onto $F$, and $(Tf)^{-1}=g\,(f\circ\vartheta)^{-1}$ gives $\psi([f])=[(Tf)^{-1}]$ for all $[f]\in P(E_+^\circ)$.

In either case, $\psi([f])=[(Tf)^\varepsilon]$ for all $[f]\in P(E_+^\circ)$, where $T\colon E\to F$ is a lattice isomorphism.
\endgroup

\begingroup\relpenalty=10000 \binoppenalty=10000
$(\Leftarrow)$. Let $T\colon E\to F$ be a lattice isomorphism. If $\varepsilon=1$, then $\psi$ is a surjective Hilbert isometry by Lemma~\ref{lem:lattice-isomorphism-order-units-isometry}(iii). Suppose now that $\varepsilon=-1$ and that $\psi([f])=[(Tf)^{-1}]$ defines a map from $P(E_+^\circ)$ into $P(F_+^\circ)$. Since $T(E_+^\circ)=F_+^\circ$ by Lemma~\ref{lem:lattice-isomorphism-order-units-isometry}(i), for every $x\in F_+^\circ$ the inverse $x^{-1}$ is a positive multiple of an element of $F_+^\circ$, and hence by Corollary~\ref{lem:inverse-isometry-well-defined}, $F$ is a bounded inverse closed subalgebra of $C(K_F)$. Hence $\psi$, being the composition of $[f]\mapsto[Tf]$ with the inversion $[x]\mapsto[x^{-1}]$ on $P(F_+^\circ)$, is a surjective Hilbert isometry by Lemma~\ref{lem:lattice-isomorphism-order-units-isometry}(iii) and Lemma~\ref{lem:inversion-order-reversing}.
\endgroup
\end{proof}

In particular, Hilbert's metric determines the vector lattice up to a lattice isomorphism.

\begin{corollary}\label{C: d_H isometric iff lattice isom}
    Let $E$ and $F$ be two Archimedean vector lattices with strong order units. Then the projective cones $P(E_+^\circ)$ and $P(F_+^\circ)$ are isometric for Hilbert's metric if and only if $E$ and $F$ are lattice isomorphic.
\end{corollary}

\begin{proof}
The forward implication follows from Theorem~\ref{T:d_H-isom between E and F}, and the reverse implication from Lemma~\ref{lem:lattice-isomorphism-order-units-isometry}(iii).
\end{proof}

\subsection{The isometry group}\label{subsec:isometry-group}

For an Archimedean vector lattice $E$, we denote the group of lattice automorphisms of $E$ by $\aut(E)$. If $E$ has a strong order unit, we denote by $\isom(P(E_+^\circ))$ the group of surjective Hilbert isometries of $P(E_+^\circ)$, and by $[\aut(E)]$ the group of all maps
\[
[T]\colon P(E_+^\circ)\to P(E_+^\circ),\qquad [f]\mapsto[Tf],
\]
with $T\in\aut(E)$. By Lemma~\ref{lem:lattice-isomorphism-order-units-isometry}(iii), we have $[\aut(E)]\subseteq\isom(P(E_+^\circ))$.
 
The following consequence of Theorem~\ref{T:d_H-isom between E and F} shows that the case $\varepsilon=-1$ in that theorem can only occur when $F$ is a bounded inverse closed $f$-algebra.

\begin{corollary}\label{cor:Aut-Isom-agree}
Let $E$ and $F$ be Archimedean vector lattices with strong order units, and suppose that $F$ is not a bounded inverse closed $f$-algebra. Then a map $\psi\colon P(E_+^\circ)\to P(F_+^\circ)$ is a surjective Hilbert isometry if and only if there exists a lattice isomorphism $T\colon E\to F$ such that $\psi([f])=[Tf]$ for all $[f]\in P(E_+^\circ)$. In particular, if $E$ is not a bounded inverse closed $f$-algebra, then $\isom(P(E_+^\circ))=[\aut(E)]$.
\end{corollary}

\begin{proof}
The reverse implication follows from Lemma~\ref{lem:lattice-isomorphism-order-units-isometry}(iii). Conversely, by Theorem~\ref{T:d_H-isom between E and F} there exist $\varepsilon\in\{\pm1\}$ and a lattice isomorphism $T\colon E\to F$ such that $\psi([f])=[(Tf)^\varepsilon]$ for all $[f]\in P(E_+^\circ)$. If $\varepsilon=-1$, then, as in the proof of the reverse implication of Theorem~\ref{T:d_H-isom between E and F}, we have $x^{-1}\in F$ for every $x\in F_+^\circ$, so that $F$ is a bounded inverse closed $f$-algebra by Corollary~\ref{lem:inverse-isometry-well-defined}, which contradicts the assumption. Hence $\varepsilon=1$.
\end{proof}

If $E$ is a bounded inverse closed $f$-algebra, then the inversion is an additional surjective Hilbert isometry, and the following theorem shows that it lies outside $[\aut(E)]$ as soon as $\dim E\geq3$.

\begingroup\relpenalty=10000 \binoppenalty=10000
\begin{theorem}\label{thm:isom-group-isomorphism}
Let $E$ be an Archimedean vector lattice with strong order unit, and let $\iota([f]):=[f^{-1}]$ denote the inversion on $P(E_+^\circ)$ whenever $E$ is a bounded inverse closed $f$-algebra. If $\dim E\geq3$ and $E$ is a bounded inverse closed $f$-algebra, then $[\aut(E)]$ is a normal subgroup of index $2$ in $\isom(P(E_+^\circ))$, and
\[
\isom(P(E_+^\circ))=[\aut(E)]\rtimes\langle\iota\rangle\cong[\aut(E)]\rtimes\Z/2\Z.
\]
Otherwise, $\isom(P(E_+^\circ))=[\aut(E)]$.
\end{theorem}

\begin{proof}
We regard $E$ as a norm dense sublattice of $C(K)$ containing $\mathbb 1$. If $E$ is not a bounded inverse closed $f$-algebra, then $\isom(P(E_+^\circ))=[\aut(E)]$ by Corollary~\ref{cor:Aut-Isom-agree}. Suppose now that $E$ is a bounded inverse closed $f$-algebra. By Corollary~\ref{lem:inverse-isometry-well-defined} and Lemma~\ref{lem:inversion-order-reversing}, $\iota$ is a well-defined surjective Hilbert isometry of $P(E_+^\circ)$ with $\iota^2=\mathrm{id}$. We claim that
\[
\isom(P(E_+^\circ))=[\aut(E)]\cup\iota\,[\aut(E)].
\]
Indeed, let $\psi\in\isom(P(E_+^\circ))$. By Theorem~\ref{T:d_H-isom between E and F}, there is $T\in\aut(E)$ such that either $\psi=[T]$, or $\psi([f])=[(Tf)^{-1}]$ for all $[f]\in P(E_+^\circ)$; in the latter case, $\psi=\iota\circ[T]$. Conversely, $\iota$ and every $[T]$ with $T\in\aut(E)$ belong to the group $\isom(P(E_+^\circ))$, and hence so does $\iota\circ[T]$. In particular, the right-hand side is a group, since it equals $\isom(P(E_+^\circ))$.

If $\dim E=1$, then $P(E_+^\circ)$ consists of a single point, so $\isom(P(E_+^\circ))=[\aut(E)]$ is trivial. If $\dim E=2$, then $E$ is lattice isomorphic to $\R^2$ with the coordinatewise order by \cite[Theorem~26.11]{LZ71}. Since $(x^{-1},y^{-1})=(xy)^{-1}(y,x)$ for all $(x,y)\in(\R^2)_+^\circ$, the inversion $\iota$ is induced by the lattice automorphism that interchanges the coordinates. Hence $\iota\in[\aut(E)]$, and so $\isom(P(E_+^\circ))=[\aut(E)]$.

Finally, suppose that $\dim E\geq3$. We claim that $\iota\notin[\aut(E)]$. Suppose, to the contrary, that there is $T\in\aut(E)$ with $[Tf]=[f^{-1}]$ for all $f\in E_+^\circ$. Then for every $f\in E_+^\circ$ there is $\alpha_f>0$ with $Tf=\alpha_f f^{-1}$; in particular, $T\mathbb 1=\alpha_{\mathbb 1}\mathbb 1$. Since $\dim E\geq3$, the space $K$ contains three distinct points $t_1,t_2,t_3$. By Urysohn's lemma, for each $i=2,3$ there is $g_i\in C(K)$ satisfying $g_i(t_j)=\delta_{ij}$ for $j=1,2,3$. Since $E$ is dense in $C(K)$, we may choose $h_i\in E$ such that $\|g_i-h_i\|_\infty<1/3$. Consequently, the functions
\[
f_i:=\frac{(h_i-\frac13\mathbb 1)^+}{h_i(t_i)-\frac13}\in E_+,
\qquad i=2,3,
\]
satisfy $f_i(t_j)=\delta_{ij}$ for $j=1,2,3$. Therefore $f:=\mathbb 1+f_2+2f_3$ satisfies $f(t_i)=i$ for $i=1,2,3$. Since $f\geq\mathbb 1$, we have $f\in E_+^\circ$.

By linearity of $T$,
\[
\alpha_{f+\mathbb 1}(f+\mathbb 1)^{-1}=T(f+\mathbb 1)=Tf+T\mathbb 1=\alpha_f f^{-1}+\alpha_{\mathbb 1}\mathbb 1.
\]
Evaluating at $t_i$ and dividing by $\alpha_{f+\mathbb 1}$ gives $\frac{1}{i+1}=\frac{\alpha}{i}+\beta$ for $i=1,2,3$, where $\alpha:=\alpha_f/\alpha_{f+\mathbb 1}$ and $\beta:=\alpha_{\mathbb 1}/\alpha_{f+\mathbb 1}$. The equations for $i=1,2$ give $\alpha=\frac13$ and $\beta=\frac16$, but then $\frac14 = \frac{\alpha}{3}+\beta=\frac{5}{18}$, which is absurd. This proves the claim.

Since $\iota\notin[\aut(E)]$, the cosets $[\aut(E)]$ and $\iota\,[\aut(E)]$ are distinct, so $[\aut(E)]$ has index $2$ in $\isom(P(E_+^\circ))$ and is therefore a normal subgroup. Moreover, $\iota\neq\mathrm{id}$ because $\mathrm{id}\in[\aut(E)]$, so $\iota^2=\mathrm{id}$ implies that $\langle\iota\rangle=\{\mathrm{id},\iota\}$ has order $2$, and $\langle\iota\rangle\cap[\aut(E)]=\{\mathrm{id}\}$. Since $\isom(P(E_+^\circ))=[\aut(E)]\,\langle\iota\rangle$, it follows by \cite[Section~5.5, Theorem~12]{DumFoot04} that
\[
\isom(P(E_+^\circ))=[\aut(E)]\rtimes\langle\iota\rangle\cong[\aut(E)]\rtimes\Z/2\Z.\qedhere
\]
\end{proof}
\endgroup

\section{Applications to function lattices}\label{sec:applications}

In this section we apply the results of Section~\ref{sec:main} to concrete norm dense sublattices of spaces of continuous functions, namely to H\"older continuous functions, absolutely continuous functions, and continuous functions of bounded variation. For each of these spaces, we describe the surjective isometries for Hilbert's metric on the projective cone. By Theorem~\ref{T:d_H-isom between E and F}, this amounts to describing the lattice isomorphisms between these spaces. For dense sublattices of $C(K)$ that are bounded inverse closed $f$-subalgebras, Theorem~\ref{thm:hilbert-isometries-dense-bic} below reduces this further to determining the homeomorphisms $\vartheta$ whose composition operators $C_\vartheta f:=f\circ\vartheta$ map one of these spaces onto the other. All three spaces are bounded inverse closed $f$-subalgebras, so that Theorem~\ref{thm:hilbert-isometries-dense-bic} applies to each of them. In Subsection~\ref{subsec:comparison} we compare the projective cones of these spaces as metric spaces. We show that the projective cones of these spaces of functions on a compact interval all have the same metric completion but are not isometric, and that Hilbert isometries of these projective cones are uniquely determined by their values on three specifically chosen rays.

Throughout this section, $K$ and $L$ are compact Hausdorff spaces. For a homeomorphism $\vartheta\colon L\to K$, the composition operator $C_\vartheta\colon C(K)\to C(L)$, $C_\vartheta f:=f\circ\vartheta$, is a lattice isomorphism with inverse $C_{\vartheta^{-1}}$, since $f\circ\vartheta\leq g\circ\vartheta$ holds if and only if $f\leq g$. Consequently, if $E\subseteq C(K)$ and $F\subseteq C(L)$ are vector sublattices with $C_\vartheta(E)=F$, then $C_\vartheta$ restricts to a lattice isomorphism from $E$ onto $F$ with inverse $C_{\vartheta^{-1}}$. If, in addition, $E$ and $F$ contain the constant functions, then $E$ is a subalgebra of $C(K)$ if and only if $F$ is a subalgebra of $C(L)$, and $E$ is bounded inverse closed if and only if $F$ is, because 
\begin{equation}\label{eq:multiplicative composition operator}
    (f_1\circ\vartheta)(f_2\circ\vartheta)=(f_1f_2)\circ\vartheta, \quad (f\circ\vartheta)^{-1}=f^{-1}\circ\vartheta,
\end{equation} and $f\geq\mathbb 1_K$ if and only if $f\circ\vartheta\geq\mathbb 1_L$. In Subsections~\ref{subsec:ac} through \ref{subsec:comparison}, $K$ and $L$ are compact intervals, and we use repeatedly that every homeomorphism between compact intervals is strictly monotone, as a consequence of the intermediate value theorem.

The following theorem reduces the description of the surjective Hilbert isometries between projective cones of norm dense bounded inverse closed $f$-subalgebras to the determining what the homeomorphisms $\vartheta$ are for which $C_\vartheta$ maps one of the spaces onto the other. Its proof relies on \cite[Theorem~2.34]{AB06}, the classical description of the lattice isomorphisms between spaces of continuous functions on compact Hausdorff spaces as weighted composition operators.

\begin{theorem}\label{thm:hilbert-isometries-dense-bic}
Let $E\subseteq C(K)$ and $F\subseteq C(L)$ be norm dense bounded inverse closed $f$-subalgebras containing the constant one function.
\begin{enumerate}
\item A linear map $T\colon E\to F$ is a lattice isomorphism if and only if there exist a function $g\in F_+^\circ$ and a homeomorphism $\vartheta\colon L\to K$ with $C_\vartheta(E)=F$ such that $Tf=g\,(f\circ\vartheta)$ for all $f\in E$.
\item A map $\psi\colon P(E_+^\circ)\to P(F_+^\circ)$ is a surjective Hilbert isometry if and only if there exist $\varepsilon\in\{\pm1\}$, a homeomorphism $\vartheta\colon L\to K$ with $C_\vartheta(E)=F$, and a function $g\in F_+^\circ$ such that
\[
\psi([f])=[(g\,(f\circ\vartheta))^\varepsilon]\qquad\text{for all }[f]\in P(E_+^\circ).
\]
\item $P(E_+^\circ)$ and $P(F_+^\circ)$ are isometric for Hilbert's metric if and only if there is a homeomorphism $\vartheta\colon L\to K$ with $C_\vartheta(E)=F$.
\end{enumerate}
\end{theorem}

\begin{proof}
(i). Suppose that $T$ is a lattice isomorphism. Since we have $|Tf|\leq\|f\|_\infty\,T\mathbb 1_K$ for all $f\in E$, the map $T$ is bounded, and so it extends uniquely to a bounded linear operator $T\colon C(K)\to C(L)$. Since the lattice operations are continuous with respect to the supremum norm, this extension is again a lattice homomorphism, and it is bijective, since the extension of $T^{-1}$, obtained in the same way, is its inverse. Hence it is a lattice isomorphism. By \cite[Theorem~2.34]{AB06}, there exist a function $g\in C(L)_+^\circ$ and a homeomorphism $\vartheta\colon L\to K$ such that
\begin{equation}\label{eq:Banach-Stone-representation-T}
Tf=g\,(f\circ\vartheta)\qquad\text{for all }f\in C(K).
\end{equation}
Since $\mathbb 1_K\in E$, we have $g=T\mathbb 1_K\in F_+^\circ$ by Lemma~\ref{lem:lattice-isomorphism-order-units-isometry}(i), and hence $g^{-1}\in F$ by Corollary~\ref{lem:inverse-isometry-well-defined}. Consequently, $f\circ\vartheta=g^{-1}\,Tf\in F$ for every $f\in E$, that is, $C_\vartheta(E)\subseteq F$.

For the reverse inclusion, note that \eqref{eq:Banach-Stone-representation-T} gives $T^{-1}q=(g\circ\vartheta^{-1})^{-1}(q\circ\vartheta^{-1})$ for all $q\in C(L)$. Hence the lattice isomorphism $T^{-1}\colon F\to E$ has the same form, with the homeomorphism $\vartheta^{-1}$, and the argument above shows that $C_{\vartheta^{-1}}(F)\subseteq E$. Applying $C_\vartheta$ yields $F\subseteq C_\vartheta(E)$, and therefore $C_\vartheta(E)=F$.

Conversely, let $\vartheta$ and $g$ be as in the statement. As observed above, the map $C_\vartheta$ is a lattice isomorphism of $E$ onto $F$. By Corollary~\ref{lem:inverse-isometry-well-defined}, we have $g^{-1}\in F$, so multiplication by $g$ is a lattice automorphism of $F$. Hence $f\mapsto g\,(f\circ\vartheta)$ is a lattice isomorphism from $E$ onto $F$.

(ii). Since $E$ and $F$ are bounded inverse closed $f$-algebras, Theorem~\ref{T:d_H-isom between E and F} shows that the surjective Hilbert isometries from $P(E_+^\circ)$ onto $P(F_+^\circ)$ are exactly the maps $[f]\mapsto[(Tf)^\varepsilon]$ with $\varepsilon\in\{\pm1\}$ and $T\colon E\to F$ a lattice isomorphism. The claim therefore follows from (i). Since $F_+^\circ$ is closed under inversion by Corollary~\ref{lem:inverse-isometry-well-defined}, these maps can equivalently be written as $[f]\mapsto[g\,(f\circ\vartheta)^\varepsilon]$ with $g\in F_+^\circ$, as in Theorem~\ref{thm:Banach-Stone-projective}.

(iii). This follows from (ii), since the maps in (ii) exist if and only if there is a homeomorphism $\vartheta\colon L\to K$ with $C_\vartheta(E)=F$.
\end{proof}

Hence, for each of the spaces considered in Subsections~\ref{subsec:holder} through \ref{subsec:cbv}, it only remains to determine the homeomorphisms $\vartheta$ for which $C_\vartheta(E)=F$. We emphasize that this equivalently is a question about unital algebra isomorphisms. Indeed, if $T\colon E\to F$ is a lattice isomorphism with $Tf=g\,(f\circ\vartheta)$ as in Theorem~\ref{thm:hilbert-isometries-dense-bic}(i), then $g^{-1}\in F$, since $F$ is bounded inverse closed, and $C_\vartheta f = g^{-1}Tf$. It maps the unit $\mathbb 1_K$ to the unit $\mathbb 1_L$, and is multiplicative by the identities  \eqref{eq:multiplicative composition operator}. Hence $C_\vartheta$ is an algebra isomorphism. Conversely, every algebra isomorphism is positive. Indeed, note that $f \in E_+$ if and only if $f + \lambda\mathbb 1_K$ is invertible for all $\lambda >0$, which is a property preserved by an algebra isomorphism. By \cite[Corollary~5.5]{HdP84} this algebra isomorphism is a unital lattice isomorphism, which brings us back to the composition operator $C_\vartheta$ by Theorem~\ref{thm:hilbert-isometries-dense-bic}(i). Thus, up to a weight, the lattice isomorphisms between such spaces are precisely the unital algebra isomorphisms.

For the spaces considered below, these are known. For H\"older continuous functions the relevant characterization of $\vartheta$ is due to Garrido and Jaramillo \cite{GJ04}, for absolutely continuous functions to Al-shakarchi and Doust \cite{AlShakarchiDoust2021}, and for continuous functions of bounded variation it is a straightforward verification. The contribution of this section lies in the resulting descriptions of the surjective Hilbert isometries in Theorems~\ref{thm:hilbert-isometry-Holder-spaces}, \ref{thm:hilbert-isometries-absolutely-continuous}, and~\ref{thm:hilbert-isometries-bounded-variation}, in which the weights and the inversion enter through Theorem~\ref{T:d_H-isom between E and F}, in the examples showing that these descriptions differ for the three classes of functions, and in the comparison of the resulting projective cones in Subsection~\ref{subsec:comparison}.

\subsection{H\"{o}lder continuous functions}\label{subsec:holder}

Let $(K,d)$ be a compact metric space and let $0<\alpha\leq1$. We denote by $C^{0,\alpha}(K)$ the set of all functions $f\colon K\to\mathbb R$ such that
\begin{equation*}\label{eq:Holder-continuous-function}
L_{\alpha}(f):=\sup_{x\neq y} \frac{|f(x)-f(y)|}{d(x,y)^\alpha}<\infty.
\end{equation*}
Such functions are called \emph{$\alpha$-H\"older continuous}, and for $\alpha=1$ they are the Lipschitz functions on $K$. Since $(s+t)^\alpha\leq s^\alpha+t^\alpha$ for all $s,t\geq0$, the function $d^\alpha$ is again a metric on $K$, which induces the same topology as $d$, and $C^{0,\alpha}(K)$ is precisely the space of Lipschitz functions on $(K,d^\alpha)$. The surjective linear isometries of $C^{0,\alpha}(K)$ for the norm $\max\{\|f\|_\infty,L_\alpha(f)\}$ were described by Jarosz and Pathak \cite[Example~3]{JP88}, who also treated the norm $\|f\|_\infty+L_\alpha(f)$ for complex-valued functions \cite[Example~8]{JP88}; see also the references therein. In both cases, every surjective linear isometry is a weighted composition operator $f\mapsto\chi\cdot(f\circ\vartheta)$, where $\vartheta$ is a homeomorphism and $\chi$ is a continuous function with $|\chi|=1$. For real-valued functions, $\chi$ takes only the values $\pm1$ and is constant on each connected component of $K$. By contrast, the surjective Hilbert isometries described in Theorem~\ref{thm:hilbert-isometry-Holder-spaces} below also involve arbitrary strictly positive weights and the inversion $[f]\mapsto[f^{-1}]$.

The space $C^{0,\alpha}(K)$ is a norm dense bounded inverse closed $f$-subalgebra of $C(K)$ that contains $\mathbb 1_K$. Indeed, for $f,g\in C^{0,\alpha}(K)$ and $x,y\in K$ we have
\[
|(fg)(x)-(fg)(y)|\leq\bigl(\|f\|_\infty L_\alpha(g)+\|g\|_\infty L_\alpha(f)\bigr)\,d(x,y)^\alpha
\]
and $\bigl||f(x)|-|f(y)|\bigr|\leq|f(x)-f(y)|$, so $C^{0,\alpha}(K)$ is a subalgebra and a vector sublattice of $C(K)$ containing the constant functions. If $f\geq\mathbb 1_K$, then
\[
\bigl|f(x)^{-1}-f(y)^{-1}\bigr|=\frac{|f(x)-f(y)|}{f(x)f(y)}\leq L_\alpha(f)\,d(x,y)^\alpha,
\]
so $f^{-1}\in C^{0,\alpha}(K)$. Finally, for $z\in K$ the function $x\mapsto d(z,x)^\alpha$ belongs to $C^{0,\alpha}(K)$ by the triangle inequality for $d^\alpha$, so $C^{0,\alpha}(K)$ separates the points of $K$ and is norm dense in $C(K)$ by the lattice version of the Stone--Weierstrass theorem \cite[Theorem~2.1.1]{Mey91}.

The composition operators between spaces of H\"older continuous functions are well understood. Garrido and Jaramillo \cite[Theorem~3.10]{GJ04} showed that for complete metric spaces $X$ and $Y$, every unital vector lattice isomorphism between the spaces of Lipschitz functions on $Y$ and on $X$ is of the form $f\mapsto f\circ\vartheta$ for a bi-Lipschitz homeomorphism $\vartheta\colon X\to Y$; see also \cite[Theorem~5.1]{She63} for the corresponding result for algebra isomorphisms. For compact metric spaces $(K,d_K)$ and $(L,d_L)$ and a homeomorphism $\vartheta\colon L\to K$, this yields that
\begin{equation}\label{eq:holder-composition-bilipschitz}
C_\vartheta(C^{0,\alpha}(K))=C^{0,\alpha}(L)\qquad\text{if and only if $\vartheta$ is bi-Lipschitz.}
\end{equation}
Indeed, if $\vartheta$ is bi-Lipschitz, then $L_\alpha(f\circ\vartheta)\leq c^\alpha L_\alpha(f)$ for every $f\in C^{0,\alpha}(K)$, where $c$ is a Lipschitz constant of $\vartheta$; the same holds for $\vartheta^{-1}$. Conversely, suppose that $C_\vartheta(C^{0,\alpha}(K))=C^{0,\alpha}(L)$. Then $C_\vartheta$ is a unital lattice isomorphism between the spaces of Lipschitz functions on $(K,d_K^\alpha)$ and on $(L,d_L^\alpha)$. By \cite[Theorem~3.10]{GJ04}, there is a homeomorphism $\vartheta'\colon L\to K$, which is bi-Lipschitz with respect to $d_L^\alpha$ and $d_K^\alpha$, such that $f\circ\vartheta=f\circ \vartheta'$ for all $f\in C^{0,\alpha}(K)$. Since $C^{0,\alpha}(K)$ separates the points of $K$, we have $\vartheta'=\vartheta$, and $\vartheta$ is bi-Lipschitz for $d_K^\alpha$ and $d_L^\alpha$ if and only if it is bi-Lipschitz for $d_K$ and $d_L$.

We can now describe the surjective Hilbert isometries between projective cones of H\"older continuous functions.

\begin{theorem}\label{thm:hilbert-isometry-Holder-spaces}
Let $(K,d_K)$ and $(L,d_L)$ be compact metric spaces and let $0<\alpha\leq 1$. A map $\psi\colon P(C^{0,\alpha}(K)_+^\circ)\to P(C^{0,\alpha}(L)_+^\circ)$ is a surjective Hilbert isometry if and only if there exist $\varepsilon\in\{\pm1\}$, a function $g\in C^{0,\alpha}(L)_+^\circ$, and a bi-Lipschitz homeomorphism $\vartheta\colon L\to K$ such that
\[
\psi([f])=[(g\,(f\circ\vartheta))^\varepsilon]\qquad\text{for all }[f]\in P(C^{0,\alpha}(K)_+^\circ).
\]
In particular, the projective cones $P(C^{0,\alpha}(K)_+^\circ)$ and $P(C^{0,\alpha}(L)_+^\circ)$ are isometric for Hilbert's metric if and only if $K$ and $L$ are bi-Lipschitz homeomorphic.
\end{theorem}

\begin{proof}
Put $E:=C^{0,\alpha}(K)$ and $F:=C^{0,\alpha}(L)$. As observed above, these spaces satisfy the hypotheses of Theorem~\ref{thm:hilbert-isometries-dense-bic}, and by \eqref{eq:holder-composition-bilipschitz} a homeomorphism $\vartheta\colon L\to K$ satisfies $C_\vartheta(E)=F$ if and only if it is bi-Lipschitz. Both assertions therefore follow from Theorem~\ref{thm:hilbert-isometries-dense-bic}(ii) and (iii).
\end{proof}

\begin{example}\label{ex:homeomorphic-not-bilipschitz}
{By Theorem~\ref{thm:Banach-Stone-projective}, the projective cones $P(C(K)_+^\circ)$ and $P(C(L)_+^\circ)$ are isometric if and only if $K$ and $L$ are homeomorphic. Theorem~\ref{thm:hilbert-isometry-Holder-spaces} shows that for H\"older continuous functions the metric structure of $K$ and $L$ matters as well, since homeomorphic compact metric spaces need not be bi-Lipschitz homeomorphic. For instance, let $0<\beta<1$, and let $K=[0,1]$ with the usual metric $d(x,y)=|x-y|$ and $L=[0,1]$ with the metric $d^\beta$. Suppose that $\vartheta\colon L\to K$ is a bi-Lipschitz homeomorphism. Since $d^\beta$ induces the same topology as $d$ on $[0,1]$, the map $\vartheta$ is a homeomorphism of $[0,1]$ for $d$, and is therefore strictly monotone with $\{\vartheta(0),\vartheta(1)\}=\{0,1\}$. Moreover, there is $c>0$ with $|\vartheta(x)-\vartheta(y)|\geq c\,|x-y|^\beta$ for all $x,y\in[0,1]$. For $n\in\mathbb N$ and $x_i:=i/n$ for $i = 1, \dots, n$, this yields
\[
1=|\vartheta(1)-\vartheta(0)|=\sum_{i=1}^n|\vartheta(x_i)-\vartheta(x_{i-1})|\geq c\,n^{1-\beta},
\]
which is impossible for large enough $n$. Hence $K$ and $L$ are not bi-Lipschitz homeomorphic, and Theorem~\ref{thm:hilbert-isometry-Holder-spaces} shows that $P(C^{0,\gamma}(K)_+^\circ)$ and $P(C^{0,\gamma}(L)_+^\circ)$ are not isometric for Hilbert’s metric for any $0<\gamma\le1$. Now, a function on $[0,1]$ is $\gamma$-H\"{o}lder continuous for the metric $d^\beta$ if and only if it is $\beta\gamma$-H\"{o}lder continuous for the metric $d$, so that $C^{0,\gamma}(L)=C^{0,\beta\gamma}(K)$. Since every pair of exponents $0<\alpha<\gamma\le1$ is of the form $\alpha=\beta\gamma$ with $\beta=\alpha/\gamma\in(0,1)$, it follows that, for $[0,1]$ with the metric $d$, the projective cones $P(C^{0,\alpha}[0,1]_+^\circ)$ and $P(C^{0,\gamma}[0,1]_+^\circ)$ are not isometric for Hilbert’s metric whenever $0<\alpha<\gamma\le1$. Thus the projective cone of the space of H\"{o}lder continuous functions on $[0,1]$ determines the H\"{o}lder exponent.}
\end{example}

\subsection{Absolutely continuous functions}\label{subsec:ac}

Let $a<b$ be real numbers. We denote by $\ac[a,b]$ the set of all functions $f\colon[a,b]\to\mathbb R$ with the property that for each $\varepsilon>0$ there exists $\delta>0$ such that
\[
\sum_{i=1}^n |f(t_i)-f(s_i)| < \varepsilon
\]
for every finite family of pairwise disjoint intervals $(s_1,t_1),\ldots,(s_n,t_n)$ in $[a,b]$ with $\sum_{i=1}^n(t_i-s_i)<\delta$. Functions in $\ac[a,b]$ are called \emph{absolutely continuous}. For a function $f\colon[a,b]\to\mathbb R$, its \emph{total variation} is
\[
V_a^b(f):=\sup\Bigl\{\sum_{i=1}^n|f(x_i)-f(x_{i-1})|:\ n\in\mathbb N,\ a=x_0<x_1<\cdots<x_n=b\Bigr\}.
\]
If $V_a^b(f)<\infty$, then $f$ is said to be of \emph{bounded variation}. Every absolutely continuous function is of bounded variation. The surjective linear isometries of $\ac[a,b]$ were described by Jarosz and Pathak \cite[Example~2]{JP88} for the norm $\max\{\|f\|_\infty,V_a^b(f)\}$, and by Pathak \cite[Theorem~2.10]{Pat82} for the norm $\|f\|_\infty+V_a^b(f)$ in the case of complex-valued functions. In both cases, every surjective linear isometry between such spaces is of the form $f\mapsto\sigma\,(f\circ\vartheta)$, where $\sigma$ is a scalar with $|\sigma|=1$, so that $\sigma=\pm1$ for real-valued functions, and $\vartheta$ is a homeomorphism between the underlying intervals. Moreover, Dubarbie \cite[Theorem~3.6]{Dub10} showed that every disjointness preserving linear bijection between spaces of absolutely continuous functions, and in particular every lattice isomorphism, is a weighted composition operator induced by a homeomorphism, and Al-shakarchi and Doust \cite[Theorem~4.7]{AlShakarchiDoust2021} characterized the homeomorphisms $\vartheta$ for which the composition operator $C_\vartheta$ is an algebra isomorphism between spaces of absolutely continuous functions on compact subsets of the plane. For intervals, their characterization is Lemma~\ref{lem:composition-operator-absolutely-continuous} below. By contrast, the surjective Hilbert isometries described in Theorem~\ref{thm:hilbert-isometries-absolutely-continuous} below also involve arbitrary strictly positive weights and the inversion $[f]\mapsto[f^{-1}]$.

The space $\ac[a,b]$ is a norm dense bounded inverse closed $f$-subalgebra of $C[a,b]$ that contains $\mathbb 1$. Indeed, absolutely continuous functions are bounded, and for all $f,g\in\ac[a,b]$ and all $x,y\in[a,b]$ we have
\[
|(fg)(x)-(fg)(y)|\leq\|f\|_\infty|g(x)-g(y)|+\|g\|_\infty|f(x)-f(y)|
\]
and $\bigl||f(x)|-|f(y)|\bigr|\leq|f(x)-f(y)|$, so $\ac[a,b]$ is a subalgebra and a vector sublattice of $C[a,b]$ containing the constant functions. If $f\geq\mathbb 1$, then
\[
|f(x)^{-1}-f(y)^{-1}|\leq|f(x)-f(y)|,
\]
so $f^{-1}\in\ac[a,b]$. Finally, $\ac[a,b]$ contains the Lipschitz functions on $[a,b]$, which are norm dense in $C[a,b]$ by Subsection~\ref{subsec:holder}.

\begin{lemma}[cf.\ {\cite[Theorem~4.7]{AlShakarchiDoust2021}}]\label{lem:composition-operator-absolutely-continuous}
Let $a<b$ and $c<d$ be real numbers, and let $\vartheta\colon[c,d]\to[a,b]$ be a homeomorphism. Then $C_\vartheta(\ac[a,b])=\ac[c,d]$ if and only if both $\vartheta$ and $\vartheta^{-1}$ are absolutely continuous.
\end{lemma}

Lemma~\ref{lem:composition-operator-absolutely-continuous} is the interval case of \cite[Theorem~4.7]{AlShakarchiDoust2021}, which is formulated for the algebras of absolutely continuous functions and of functions of bounded variation on compact subsets $S$ of the plane, see \cite[Section~2]{AlShakarchiDoust2021}. These are Banach algebras for the norm $\|f\|_\infty+\operatorname{var}(f,S)$, where $\operatorname{var}(f,S)$ denotes their planar notion of variation. For a compact interval $[a,b]$, these spaces coincide with the classical ones and $\operatorname{var}(f,[a,b])=V_a^b(f)$, see \cite[Section~2]{AlShakarchiDoust2021}; every homeomorphism between compact intervals preserves the total variation, being monotone, and therefore induces an algebra isomorphism between the spaces of all functions of bounded variation, and algebra isomorphisms between these spaces are automatically continuous. The implication from right to left also follows from the elementary observation that, for an absolutely continuous homeomorphism $\vartheta$, the composition $f\circ\vartheta$ is absolutely continuous for every absolutely continuous $f$, as $\vartheta$ maps disjoint intervals of small total length to disjoint intervals of small total length.

We can now describe the surjective Hilbert isometries between projective cones of absolutely continuous functions.

\begin{theorem}\label{thm:hilbert-isometries-absolutely-continuous}
Let $a<b$ and $c<d$ be real numbers, and let $\psi$ be a map from $P(\ac[a,b]_+^\circ)$ to $P(\ac[c,d]_+^\circ)$. Then $\psi$ is a surjective Hilbert isometry if and only if there exist $\varepsilon\in\{\pm 1\}$, a function $g\in\ac[c,d]_+^\circ$, and a homeomorphism $\vartheta\colon[c,d]\to[a,b]$ such that both $\vartheta$ and $\vartheta^{-1}$ are absolutely continuous and
\[
\psi([f])=[(g\,(f\circ\vartheta))^\varepsilon]\qquad\text{for all }[f]\in P(\ac[a,b]_+^\circ).
\]
In particular, $P(\ac[a,b]_+^\circ)$ and $P(\ac[c,d]_+^\circ)$ are isometric for Hilbert's metric.
\end{theorem}

\begin{proof}
Put $E:=\ac[a,b]$ and $F:=\ac[c,d]$. As observed above, these spaces satisfy the hypotheses of Theorem~\ref{thm:hilbert-isometries-dense-bic}, and by Lemma~\ref{lem:composition-operator-absolutely-continuous} a homeomorphism $\vartheta$ from $[c,d]$ onto $[a,b]$ satisfies $C_\vartheta(E)=F$ if and only if $\vartheta$ and $\vartheta^{-1}$ are absolutely continuous. The first assertion therefore follows from Theorem~\ref{thm:hilbert-isometries-dense-bic}(ii). For the second assertion, by Theorem~\ref{thm:hilbert-isometries-dense-bic}(iii) it suffices to find a homeomorphism $\vartheta$ from $[c,d]$ onto $[a,b]$ such that $\vartheta$ and $\vartheta^{-1}$ are absolutely continuous. The increasing affine map from $[c,d]$ onto $[a,b]$ satisfies this property.
\end{proof}

\begin{example}\label{ex:homeomorphism-not-ac}
By Theorem~\ref{thm:Banach-Stone-projective}, every homeomorphism $\vartheta$ of $[0,1]$ induces a surjective Hilbert isometry $[f]\mapsto[f\circ\vartheta]$ of $P(C[0,1]_+^\circ)$. For absolutely continuous functions this is no longer the case. For instance, let $\kappa\colon[0,1]\to[0,1]$ be the Cantor function, which is continuous and increasing with $\kappa(0)=0$, $\kappa(1)=1$, and almost everywhere the derivative can be taken and $\kappa'=0$. Hence $\kappa$ is not absolutely continuous; see the discussion following \cite[Proposition~1.22]{Fol99} and \cite[Section~3.5, Exercise~40]{Fol99}. Define
\[
\vartheta(x):=\tfrac12\bigl(x+\kappa(x)\bigr),\qquad x\in[0,1].
\]
Then $\vartheta$ is continuous and strictly increasing with $\vartheta(0)=0$ and $\vartheta(1)=1$, and hence a homeomorphism of $[0,1]$. Since $\kappa$ is increasing, we have $\vartheta(y)-\vartheta(x)\geq\frac12(y-x)$ for $x<y$, so that $\vartheta^{-1}$ is Lipschitz and in particular absolutely continuous. On the other hand, $\vartheta$ is not absolutely continuous, since otherwise $\kappa=2\vartheta-\mathrm{id}$ would be absolutely continuous. Consequently, the map $[f]\mapsto[f\circ\vartheta]$ does not even map $P(\ac[0,1]_+^\circ)$ into itself, as $(\mathbb 1+\mathrm{id})\circ\vartheta=\mathbb 1+\vartheta\notin\ac[0,1]$. Moreover, the absolute continuity of both $\vartheta$ and $\vartheta^{-1}$ is needed in Lemma~\ref{lem:composition-operator-absolutely-continuous}. The same homeomorphism is used in \cite[Example~4.4]{AlShakarchiDoust2021} to show that an algebra isomorphism between spaces of functions of bounded variation need not preserve absolute continuity.

A similar construction also yields an isometric embedding of
$P(\ac[0,1]_+^\circ)$ into itself for Hilbert's metric that is not surjective. Indeed, put $\eta=\vartheta^{-1}$. Since $\eta$ is Lipschitz,
composition with $\eta$ preserves absolute continuity. Moreover, since $\eta$ is a
homeomorphism, the mapping $[f]\mapsto[f\circ\eta]$ preserves Hilbert's
metric and is therefore an isometric embedding of $P(\ac[0,1]_+^\circ)$ into itself.
However, it is not surjective. Otherwise, there would exist $f\in\ac[0,1]_+^\circ$
and $\lambda>0$ such that $f\circ\eta=\lambda(\mathbb 1+\operatorname{id})$.
Composing with $\vartheta$ would give $f=\lambda(\mathbb 1+\vartheta)$,
contradicting the fact that $\vartheta$ is not absolutely continuous.
\end{example}

\subsection{Continuous functions of bounded variation}\label{subsec:cbv}

Let $a<b$ be real numbers, and recall the total variation $V_a^b(f)$ of a function $f\colon[a,b]\to\mathbb R$ defined in Subsection~\ref{subsec:ac}. We denote by $\bv[a,b]$ the set of all functions of bounded variation on $[a,b]$, and by $\cbv[a,b]:=\bv[a,b]\cap C[a,b]$ the set of all continuous functions of bounded variation on $[a,b]$. Araujo \cite{Ara99} showed that every surjective linear isometry between spaces of complex-valued functions of bounded variation on subsets of the real line, equipped with the norm $\|f\|_\infty+V(f)$, where $V(f)$ denotes the total variation, is of the form $f\mapsto\sigma\,(f\circ\vartheta)$, where $\sigma$ is a scalar with $|\sigma|=1$ and $\vartheta$ is a monotone bijection. By contrast, the surjective Hilbert isometries described in Theorem~\ref{thm:hilbert-isometries-bounded-variation} below also involve arbitrary strictly positive weights and the inversion $[f]\mapsto[f^{-1}]$.

The space $\cbv[a,b]$ is a norm dense bounded inverse closed $f$-subalgebra of $C[a,b]$ that contains $\mathbb 1$. To see this, let $f,g\in\cbv[a,b]$, and consider a partition
\[
a=x_0<x_1<\cdots<x_n=b
\]
of $[a,b]$. Since
\[
|(fg)(x_i)-(fg)(x_{i-1})|\leq\|f\|_\infty|g(x_i)-g(x_{i-1})|+\|g\|_\infty|f(x_i)-f(x_{i-1})|
\]
and $\bigl||f(x_i)|-|f(x_{i-1})|\bigr|\leq|f(x_i)-f(x_{i-1})|$ for all $i$, we obtain
\[
V_a^b(fg)\leq\|f\|_\infty V_a^b(g)+\|g\|_\infty V_a^b(f)
\qquad\text{and}\qquad
V_a^b(|f|)\leq V_a^b(f).
\]
Together with $V_a^b(f+g)\leq V_a^b(f)+V_a^b(g)$ and $V_a^b(\lambda f)=|\lambda|V_a^b(f)$, this shows that $\cbv[a,b]$ is a subalgebra and a vector sublattice of $C[a,b]$ containing the constant functions. If $f\geq\mathbb 1$, then $|f(x_i)^{-1}-f(x_{i-1})^{-1}|\leq|f(x_i)-f(x_{i-1})|$, so that $V_a^b(f^{-1})\leq V_a^b(f)$ and $f^{-1}\in\cbv[a,b]$. Finally, $\cbv[a,b]$ contains $\ac[a,b]$, which is norm dense in $C[a,b]$.

Moreover, every homeomorphism $\vartheta\colon[c,d]\to[a,b]$ satisfies
\begin{equation}\label{eq:cbv-composition}
C_\vartheta(\cbv[a,b])=\cbv[c,d].
\end{equation}
Indeed, since $\vartheta$ is strictly monotone, it maps partitions of $[c,d]$ bijectively onto partitions of $[a,b]$, up to reversing their order, so that $V_c^d(f\circ\vartheta)=V_a^b(f)$ for every function $f\colon[a,b]\to\mathbb R$, and $f\circ\vartheta$ is continuous if and only if $f$ is continuous.

We can now describe the surjective Hilbert isometries between projective cones of continuous functions of bounded variation.

\begin{theorem}\label{thm:hilbert-isometries-bounded-variation}
Let $a<b$ and $c<d$ be real numbers. A map $\psi$ from $P(\cbv[a,b]_+^\circ)$ to $P(\cbv[c,d]_+^\circ)$ is a surjective Hilbert isometry if and only if there exist $\varepsilon\in\{\pm1\}$, a function $g\in\cbv[c,d]_+^\circ$, and a homeomorphism $\vartheta\colon[c,d]\to[a,b]$ such that
\[
\psi([f])=[(g\,(f\circ\vartheta))^\varepsilon]\qquad\text{for all }[f]\in P(\cbv[a,b]_+^\circ).
\]
In particular, $P(\cbv[a,b]_+^\circ)$ and $P(\cbv[c,d]_+^\circ)$ are isometric for Hilbert's metric.
\end{theorem}

\begin{proof}
Put $E:=\cbv[a,b]$ and $F:=\cbv[c,d]$. As observed above, these spaces satisfy the hypotheses of Theorem~\ref{thm:hilbert-isometries-dense-bic}, and by \eqref{eq:cbv-composition} every homeomorphism $\vartheta$ from $[c,d]$ onto $[a,b]$ satisfies $C_\vartheta(E)=F$. Both assertions therefore follow from Theorem~\ref{thm:hilbert-isometries-dense-bic}(ii) and (iii).
\end{proof}

In contrast to the situation for absolutely continuous functions described in Example~\ref{ex:homeomorphism-not-ac}, every homeomorphism $\vartheta\colon[c,d]\to[a,b]$ induces a surjective Hilbert isometry $[f]\mapsto[f\circ\vartheta]$ from $P(\cbv[a,b]_+^\circ)$ onto $P(\cbv[c,d]_+^\circ)$, just as for spaces of continuous functions by Theorem~\ref{thm:Banach-Stone-projective}.

\subsection{Comparison of the projective cones}\label{subsec:comparison}

In this subsection we compare the projective cones of the spaces studied in Subsections~\ref{subsec:holder} through \ref{subsec:cbv} with each other and with the projective cone of $C[a,b]$. Throughout, $a<b$ and $c<d$ are real numbers, and the spaces $C^{0,\alpha}[a,b]$ of H\"older continuous functions are taken with respect to the usual metric on $[a,b]$. Since the order on each of the spaces $\ac[a,b]$, $\cbv[a,b]$, and $C[a,b]$ is the pointwise order, Hilbert's metric on each of the smaller projective cones is the restriction of Hilbert's metric on $P(C[a,b]_+^\circ)$, so that we obtain the natural isometric inclusions
\[
P(\ac[a,b]_+^\circ)\subseteq P(\cbv[a,b]_+^\circ)\subseteq P(C[a,b]_+^\circ).
\]
The same holds for the isometric inclusion of $P(C^{0,\alpha}[a,b]_+^\circ)$ into $P(C[a,b]_+^\circ)$ for every $0<\alpha\leq1$.

We start with the projective cones of H\"older continuous functions. By Example~\ref{ex:homeomorphic-not-bilipschitz}, the projective cones of H\"older continuous functions on compact intervals with different exponents are not isometric. The next result shows that they are not isometric to the projective cones of absolutely continuous functions or of continuous functions of bounded variation either, although all these projective cones have the same metric completion.

\begin{proposition}\label{prop:holder-comparison}
Let $0<\alpha\leq1$.Then the metric completion of $P(C^{0,\alpha}[a,b]_+^\circ)$ for Hilbert's metric is $P(C[a,b]_+^\circ)$, and $P(C^{0,\alpha}[a,b]_+^\circ)$ is neither isometric to $P(\ac[c,d]_+^\circ)$ nor isometric to $P(\cbv[c,d]_+^\circ)$.
\end{proposition}

\begin{proof}
Since $C^{0,\alpha}[a,b]$ is a norm dense sublattice of $C[a,b]$ containing $\mathbb 1$, it follows from Lemma~\ref{lem:dense_sublattice_dense_projectice_cone} that $P(C^{0,\alpha}[a,b]_+^\circ)$ is dense in $P(C[a,b]_+^\circ)$. Since $P(C[a,b]_+^\circ)$ is complete, as is shown in the proof of Lemma~\ref{L: unique extension of h}, we conclude that the metric completion of $P(C^{0,\alpha}[a,b]_+^\circ)$ is $P(C[a,b]_+^\circ)$. 

Let $F$ be one of the spaces $\ac[c,d]$ and $\cbv[c,d]$. By Theorem~\ref{thm:hilbert-isometries-dense-bic}(iii), it suffices to show that there is no homeomorphism $\vartheta$ from $[c,d]$ onto $[a,b]$ such that $C_\vartheta(C^{0,\alpha}[a,b])=F$. Suppose, to the contrary, that such a homeomorphism $\vartheta$ exists. Fix $t_0\in[c,d]$ and put $x_0:=\vartheta(t_0)$. The function $h(t):=|t-t_0|^{\alpha/2}$ is absolutely continuous on $[c,d]$, so that $h\in F$. Hence $h=f\circ\vartheta$ for some $f\in C^{0,\alpha}[a,b]$, that is, $f(x)=|\vartheta^{-1}(x)-t_0|^{\alpha/2}$ for all $x\in[a,b]$. Since $f(x_0)=0$, we obtain
\[
|\vartheta^{-1}(x)-\vartheta^{-1}(x_0)|^{\alpha/2}=|f(x)-f(x_0)|\leq L_\alpha(f)\,|x-x_0|^{\alpha}\qquad\text{for all }x\in[a,b],
\]
and therefore $|\vartheta^{-1}(x)-\vartheta^{-1}(x_0)|\leq L_\alpha(f)^{2/\alpha}\,|x-x_0|^2$. This observation implies that $\vartheta^{-1}$ is differentiable at $x_0$ with derivative $0$. Since $t_0\in[c,d]$ was arbitrary, the derivative of $\vartheta^{-1}$ vanishes at every point of $[a,b]$, so that $\vartheta^{-1}$ is constant, which is impossible.
\end{proof}

Although absolutely continuous functions and continuous functions of bounded variation admit different classes of Hilbert isometries, their projective cones are closely related as metric spaces. The next result shows that they also have the same metric completion, but are nevertheless not isometric.

\begin{proposition}\label{prop:ac-cbv-comparison}
The projective cones $P(\ac[a,b]_+^\circ)$ and $P(\cbv[a,b]_+^\circ)$ have the same metric completion for Hilbert's metric, namely $P(C[a,b]_+^\circ)$, whereas $P(\ac[a,b]_+^\circ)$ and $P(\cbv[c,d]_+^\circ)$ are not isometric for Hilbert's metric.
\end{proposition}

\begin{proof}
As in the proof of Proposition~\ref{prop:holder-comparison}, $P(\ac[a,b]_+^\circ)$ is dense in the complete space $P(C[a,b]_+^\circ)$, and therefore also in the intermediate space $P(\cbv[a,b]_+^\circ)$. Hence $P(C[a,b]_+^\circ)$ is the metric completion of both smaller cones.

To prove the second assertion, we apply Theorem~\ref{thm:hilbert-isometries-dense-bic}(iii) to $E:=\ac[a,b]$ and $F:=\cbv[c,d]$. It shows that $P(\ac[a,b]_+^\circ)$ and $P(\cbv[c,d]_+^\circ)$ are isometric if and only if there is a homeomorphism $\vartheta$ from $[c,d]$ onto $[a,b]$ such that
\[
C_\vartheta(\ac[a,b])=\cbv[c,d].
\]
No such homeomorphism exists. Indeed, by \eqref{eq:cbv-composition} the operator $C_\vartheta$ is a bijection from $\cbv[a,b]$ onto $\cbv[c,d]$, and $\ac[a,b]$ is a proper subspace of $\cbv[a,b]$, since the Cantor function, transferred to $[a,b]$ by an affine change of variables, is continuous and of bounded variation, but not absolutely continuous.
\end{proof}

We also note that, by Theorems~\ref{thm:hilbert-isometries-absolutely-continuous} and~\ref{thm:hilbert-isometries-bounded-variation}, every surjective Hilbert isometry between the projective cones of absolutely continuous functions is the restriction of a surjective Hilbert isometry between the projective cones of continuous functions of bounded variation, and both are restrictions of surjective Hilbert isometries between the projective cones of continuous functions, all given by the same formula. These extensions are unique by continuity, since the smaller projective cones are dense in the larger ones by Proposition~\ref{prop:ac-cbv-comparison}.

By Theorem~\ref{thm:Banach-Stone-projective}, a surjective Hilbert isometry $\psi$ between projective cones of continuous functions on intervals is determined by the weight $g$, the homeomorphism $\vartheta$, and the sign $\varepsilon$, and each of these can be read off from the image of a single projective ray. The image of $[\mathbb 1]$ is $[g]$, which determines $g$ up to a positive scalar, and this scalar does not affect $\psi$. The image of $[e^{\mathrm{id}}]$ is $[g\,e^{\varepsilon\vartheta}]$, which determines $\varepsilon\vartheta$ up to an additive constant, and for a given sign $\varepsilon$ this constant is fixed by the requirement that $\vartheta$ maps $[c,d]$ onto $[a,b]$. The remaining ambiguity is the sign $\varepsilon$, since replacing $\varepsilon$ by $-\varepsilon$ and $\vartheta$ by $a+b-\vartheta$ leaves the images of $[\mathbb 1]$ and $[e^{\mathrm{id}}]$ unchanged. This ambiguity can be resolved by considering a third ray, $[e^{\mathrm{id}^2}]$. We now make this precise. By Lemma~\ref{L: unique extension of h}, the resulting statement carries over to the projective cones of H\"{o}lder continuous functions, absolutely continuous functions, and of continuous functions of bounded variation.

\begin{proposition}\label{prop:three-ray-rigidity}
Let $a<b$ and $c<d$ be real numbers, and let $\psi_1$ and $\psi_2$ be surjective Hilbert isometries from $P(C[a,b]_+^\circ)$ onto $P(C[c,d]_+^\circ)$. If $\psi_1$ and $\psi_2$ agree on the three projective rays $[\mathbb 1]$, $[e^{\mathrm{id}}]$, and $[e^{\mathrm{id}^2}]$, then $\psi_1=\psi_2$. The same holds for surjective Hilbert isometries from $P(C^{0,\alpha}[a,b]_+^\circ)$ onto $P(C^{0,\alpha}[c,d]_+^\circ)$ for $0 < \alpha < 1$, from $P(\ac[a,b]_+^\circ)$ onto $P(\ac[c,d]_+^\circ)$ and from $P(\cbv[a,b]_+^\circ)$ onto $P(\cbv[c,d]_+^\circ)$.
\end{proposition}

\begin{proof}
By Theorem~\ref{thm:Banach-Stone-projective}, for $j=1,2$ there exist $\varepsilon_j\in\{\pm1\}$, a function $g_j\in C[c,d]_+^\circ$, and a homeomorphism $\vartheta_j\colon[c,d]\to[a,b]$ such that
\[
\psi_j([f])=[g_j(f\circ\vartheta_j)^{\varepsilon_j}]
\qquad([f]\in P(C[a,b]_+^\circ)).
\]
Since $\psi_1([\mathbb 1])=\psi_2([\mathbb 1])$, we have $[g_1]=[g_2]$. Hence there exists $\lambda>0$ such that $g_1=\lambda g_2$. Since
\[
[g_1(f\circ\vartheta_1)^{\varepsilon_1}]
=
[\lambda^{-1}g_1(f\circ\vartheta_1)^{\varepsilon_1}]
\]
for every $[f]\in P(C[a,b]_+^\circ)$, replacing $g_1$ by $\lambda^{-1}g_1$ does not change $\psi_1$. We may therefore assume that $g_1=g_2$, and we denote this function by $g$.

Since $\psi_1([e^{\mathrm{id}}])=\psi_2([e^{\mathrm{id}}])$, there exists $\nu>0$ such that
\[
g e^{\varepsilon_1\vartheta_1}
=\nu g e^{\varepsilon_2\vartheta_2}.
\]
Since $g$ is strictly positive, we may cancel it and take logarithms to obtain
\begin{equation}\label{eq:coordinate-relation}
\varepsilon_1\vartheta_1-\varepsilon_2\vartheta_2
=(\log\nu)\mathbb 1.
\end{equation}

Suppose first that $\varepsilon_1=\varepsilon_2$. Then $\vartheta_1-\vartheta_2=\mu\mathbb 1$ for some $\mu\in\mathbb R$. Since both $\vartheta_1$ and $\vartheta_2$ are bijections from $[c,d]$ onto $[a,b]$, we have
\[
[a,b]=\vartheta_1([c,d])
=\vartheta_2([c,d])+\mu
=[a+\mu,b+\mu].
\]
Hence $\mu=0$, so $\vartheta_1=\vartheta_2$ and therefore $\psi_1=\psi_2$.

Suppose instead that $\varepsilon_1=-\varepsilon_2$. Equality \eqref{eq:coordinate-relation} now implies that $\vartheta_1+\vartheta_2=\mu\mathbb 1$ for some $\mu\in\mathbb R$. Thus $\vartheta_2=\mu-\vartheta_1$. Since both functions are surjective onto $[a,b]$, we obtain
\[
[a,b]=\vartheta_2([c,d])
=\mu-\vartheta_1([c,d])
=[\mu-b,\mu-a].
\]
Hence $\mu=a+b$, so $\vartheta_2=a+b-\vartheta_1$.
Now the equality $\psi_1([e^{\mathrm{id}^2}])=\psi_2([e^{\mathrm{id}^2}])$ yields a constant $\rho>0$ such that
\[
\varepsilon_1\vartheta_1^2-\varepsilon_2\vartheta_2^2=(\log\rho)\mathbb 1.
\]
As $\varepsilon_1=-\varepsilon_2$, it follows that $\vartheta_1^2+\vartheta_2^2$ is constant. Substituting $\vartheta_2=a+b-\vartheta_1$ and using the surjectivity of $\vartheta_1$, we conclude that the function
$t\mapsto t^2+(a+b-t)^2$
is constant on $[a,b]$, which is impossible. Therefore $\varepsilon_1=\varepsilon_2$, and the first case completes the proof for projective cones of continuous functions.

Finally, let $\psi_1$ and $\psi_2$ be surjective Hilbert isometries between the projective cones of H\"{o}lder continuous functions, absolutely continuous functions, or continuous functions of bounded variation on the intervals $[a,b]$ and $[c,d]$ that agree on the three rays. Note that each of these rays belong to all of these projective cones, since $\mathbb 1$, $e^\mathrm{id}$, and $e^{\mathrm{id}^2}$ are Lipschitz. By Lemma~\ref{L: unique extension of h}, they extend to surjective Hilbert isometries $\overline\psi_1$ and $\overline\psi_2$ between the projective cones of continuous functions, which agree on the three rays as well. By the first part, $\overline\psi_1=\overline\psi_2$, and therefore $\psi_1=\psi_2$.
\end{proof}

\noindent \textbf{Acknowledgements.}
The first author was supported by the Slovenian Research and Innovation Agency program P1-0222 and grant J1-50002.\\[-1mm]

\noindent \textbf{Data Availability Statement.}
No datasets were generated or analysed during the current study.\\[-1mm]

\noindent \textbf{Conflicts of Interest.}
The authors declare that they have no conflicts of interest relevant to the content of this article.\\[-1mm]

\noindent \textbf{AI declaration.}
The authors initially used Claude Sonnet 5.5 during the development of this work to assist with identifying relevant references, examining mathematical arguments, and language editing. All mathematical arguments were independently checked and substantially revised by the authors, who take full responsibility for the content and the final version of the manuscript.

\begingroup\sloppy
\printbibliography

@book {AB06,
    AUTHOR = {Aliprantis, C. D. and Burkinshaw, O.},
     TITLE = {Positive operators},
      NOTE = {Reprint of the 1985 original},
 PUBLISHER = {Springer, Dordrecht},
      YEAR = {2006},
     PAGES = {xx+376},
      ISBN = {978-1-4020-5007-7},
   MRCLASS = {47B60 (46-01 46A40 46B42 47-01 47B65)},
  MRNUMBER = {2262133},
       DOI = {10.1007/978-1-4020-5008-4},
       URL = {https://doi.org/10.1007/978-1-4020-5008-4},
}

@article {Ara99,
    AUTHOR = {Araujo, Jes\'{u}s},
     TITLE = {Linear isometries between spaces of functions of bounded variation},
   JOURNAL = {Bull. Austral. Math. Soc.},
    VOLUME = {59},
      YEAR = {1999},
    NUMBER = {2},
     PAGES = {335--341},
       DOI = {10.1017/S0004972700032949},
}

@article {Dub10,
    AUTHOR = {Dubarbie, Luis},
     TITLE = {Separating maps between spaces of vector-valued absolutely continuous functions},
   JOURNAL = {Canad. Math. Bull.},
    VOLUME = {53},
      YEAR = {2010},
    NUMBER = {3},
     PAGES = {466--474},
       DOI = {10.4153/CMB-2010-035-7},
}

@book {DumFoot04,
    AUTHOR = {Dummit, D. S. and Foote, R. M.},
     TITLE = {Abstract algebra},
   EDITION = {3},
 PUBLISHER = {John Wiley \& Sons, Inc., Hoboken, NJ},
      YEAR = {2004},
     PAGES = {xii+932},
      ISBN = {0-471-43334-9},
   MRCLASS = {00-01 (16-01 20-01)},
  MRNUMBER = {2286236},
}

@book {Eng89,
    AUTHOR = {Engelking, Ryszard},
     TITLE = {General topology},
    SERIES = {Sigma Series in Pure Mathematics},
    VOLUME = {6},
   EDITION = {2},
 PUBLISHER = {Heldermann Verlag, Berlin},
      YEAR = {1989},
}

@book {Fol99,
    AUTHOR = {Folland, Gerald B.},
     TITLE = {Real analysis: modern techniques and their applications},
    SERIES = {Pure and Applied Mathematics (New York)},
   EDITION = {2},
 PUBLISHER = {John Wiley \& Sons, New York},
      YEAR = {1999},
}

@article {GJ04,
    AUTHOR = {Garrido, M. Isabel and Jaramillo, Jes\'{u}s A.},
     TITLE = {Homomorphisms on function lattices},
   JOURNAL = {Monatsh. Math.},
    VOLUME = {141},
      YEAR = {2004},
    NUMBER = {2},
     PAGES = {127--146},
}

@article {JP88,
    AUTHOR = {Jarosz, Krzysztof and Pathak, Vijay D.},
     TITLE = {Isometries between function spaces},
   JOURNAL = {Trans. Amer. Math. Soc.},
    VOLUME = {305},
      YEAR = {1988},
    NUMBER = {1},
     PAGES = {193--206},
}

@article {LemRoeWor18,
    AUTHOR = {Lemmens, Bas and Roelands, Mark and Wortel, Marten},
     TITLE = {Isometries of infinite dimensional {H}ilbert geometries},
   JOURNAL = {J. Topol. Anal.},
  FJOURNAL = {Journal of Topology and Analysis},
    VOLUME = {10},
      YEAR = {2018},
    NUMBER = {4},
     PAGES = {941--959},
      ISSN = {1793-5253,1793-7167},
   MRCLASS = {58B20 (22F50 46B04)},
  MRNUMBER = {3881047},
MRREVIEWER = {Leonardo\ Biliotti},
       DOI = {10.1142/S1793525318500255},
       URL = {https://doi.org/10.1142/S1793525318500255},
}

@book {LZ71,
    AUTHOR = {Luxemburg, W. A. J. and Zaanen, A. C.},
     TITLE = {Riesz spaces. {V}ol. {I}},
    SERIES = {North-Holland Mathematical Library},
 PUBLISHER = {North-Holland Publishing Co., Amsterdam-London; American
              Elsevier Publishing Co., Inc., New York},
      YEAR = {1971},
     PAGES = {xi+514},
   MRCLASS = {46A40 (06A65)},
  MRNUMBER = {511676},
MRREVIEWER = {S.\ J.\ Bernau},
}

@book {Mey91,
    AUTHOR = {Meyer-Nieberg, P.},
     TITLE = {Banach lattices},
    SERIES = {Universitext},
 PUBLISHER = {Springer-Verlag, Berlin},
      YEAR = {1991},
     PAGES = {xvi+395},
      ISBN = {3-540-54201-9},
   MRCLASS = {46B42 (46A40 47B60)},
  MRNUMBER = {1128093},
MRREVIEWER = {Yu.\ A.\ Abramovich},
       DOI = {10.1007/978-3-642-76724-1},
       URL = {https://doi.org/10.1007/978-3-642-76724-1},
}

@article {Nus88,
    AUTHOR = {Nussbaum, R. D.},
     TITLE = {Hilbert's projective metric and iterated nonlinear maps},
   JOURNAL = {Mem. Amer. Math. Soc.},
  FJOURNAL = {Memoirs of the American Mathematical Society},
    VOLUME = {75},
      YEAR = {1988},
    NUMBER = {391},
     PAGES = {iv+137},
      ISSN = {0065-9266,1947-6221},
   MRCLASS = {47H07 (47H12 58C40 92A05)},
  MRNUMBER = {961211},
MRREVIEWER = {A.\ J. B. Potter},
       DOI = {10.1090/memo/0391},
       URL = {https://doi.org/10.1090/memo/0391},
}

@article {Pat82,
    AUTHOR = {Pathak, Vijay D.},
     TITLE = {Linear isometries of spaces of absolutely continuous functions},
   JOURNAL = {Canad. J. Math.},
    VOLUME = {34},
      YEAR = {1982},
    NUMBER = {2},
     PAGES = {298--306},
}

@article {She63,
    AUTHOR = {Sherbert, Donald R.},
     TITLE = {Banach algebras of {L}ipschitz functions},
   JOURNAL = {Pacific J. Math.},
    VOLUME = {13},
      YEAR = {1963},
     PAGES = {1387--1399},
}

@article {AlShakarchiDoust2021,
    AUTHOR = {Al-shakarchi, Shaymaa and Doust, Ian},
  SORTNAME = {Alshakarchi, Shaymaa and Doust, Ian},
     TITLE = {Isomorphisms of {$BV(\sigma)$} spaces},
   JOURNAL = {Oper. Matrices},
    VOLUME = {15},
      YEAR = {2021},
    NUMBER = {4},
     PAGES = {1581--1596},
       DOI = {10.7153/oam-2021-15-101},
}

@article {Bir57,
    AUTHOR = {Birkhoff, Garrett},
     TITLE = {Extensions of {J}entzsch's theorem},
   JOURNAL = {Trans. Amer. Math. Soc.},
    VOLUME = {85},
      YEAR = {1957},
     PAGES = {219--227},
       DOI = {10.1090/S0002-9947-1957-0087058-6},
}

@article {MunozLahozTradacete2025,
    AUTHOR = {Mu{\~n}oz-Lahoz, David and Tradacete, Pedro},
     TITLE = {Banach lattice {AM}-algebras},
   JOURNAL = {Proc. Amer. Math. Soc.},
    VOLUME = {153},
      YEAR = {2025},
    NUMBER = {6},
     PAGES = {2565--2577},
}

@article {MunozLahoz2026,
    AUTHOR = {Mu{\~n}oz-Lahoz, David},
     TITLE = {{$f$}-algebra products on {AL} and {AM}-spaces},
   JOURNAL = {Math. Z.},
    VOLUME = {313},
      YEAR = {2026},
    NUMBER = {2},
     PAGES = {Paper No. 14},
       DOI = {10.1007/s00209-026-04042-3},
}

@misc {RoelandsTiersma2025,
    AUTHOR = {Roelands, Mark and Tiersma, Samuel},
     TITLE = {An order-theoretic characterization of {JB}-algebras},
      YEAR = {2025},
      NOTE = {Preprint, arXiv:2507.09526},
}

@article {HJ61,
    AUTHOR = {Henriksen, Melvin and Johnson, D. G.},
     TITLE = {On the structure of a class of archimedean lattice-ordered algebras},
   JOURNAL = {Fund. Math.},
    VOLUME = {50},
      YEAR = {1961},
     PAGES = {73--94},
       DOI = {10.4064/fm-50-1-73-94},
}

@article {HIJ61,
    AUTHOR = {Henriksen, Melvin and Isbell, J. R. and Johnson, D. G.},
     TITLE = {Residue class fields of lattice-ordered algebras},
   JOURNAL = {Fund. Math.},
    VOLUME = {50},
      YEAR = {1961},
     PAGES = {107--117},
       DOI = {10.4064/fm-50-2-107-117},
}

@article {HdP84,
    AUTHOR = {Huijsmans, C. B. and de Pagter, B.},
     TITLE = {Subalgebras and {R}iesz subspaces of an {$f$}-algebra},
   JOURNAL = {Proc. London Math. Soc. (3)},
    VOLUME = {48},
      YEAR = {1984},
    NUMBER = {1},
     PAGES = {161--174},
       DOI = {10.1112/plms/s3-48.1.161},
}

@article {Hil95,
    AUTHOR = {Hilbert, David},
     TITLE = {\"{U}ber die gerade {L}inie als k\"{u}rzeste {V}erbindung zweier {P}unkte},
   JOURNAL = {Math. Ann.},
    VOLUME = {46},
      YEAR = {1895},
     PAGES = {91--96},
}

@book {PT14,
    EDITOR = {Papadopoulos, Athanase and Troyanov, Marc},
     TITLE = {Handbook of {H}ilbert geometry},
    SERIES = {IRMA Lectures in Mathematics and Theoretical Physics},
    VOLUME = {22},
 PUBLISHER = {European Mathematical Society (EMS), Z\"{u}rich},
      YEAR = {2014},
}

@book {LN12,
    AUTHOR = {Lemmens, Bas and Nussbaum, Roger},
     TITLE = {Nonlinear {P}erron--{F}robenius theory},
    SERIES = {Cambridge Tracts in Mathematics},
    VOLUME = {189},
 PUBLISHER = {Cambridge University Press, Cambridge},
      YEAR = {2012},
}

@incollection {dlH93,
    AUTHOR = {de la Harpe, Pierre},
     TITLE = {On {H}ilbert's metric for simplices},
 BOOKTITLE = {Geometric group theory, {V}ol.\ 1 ({S}ussex, 1991)},
    SERIES = {London Math. Soc. Lecture Note Ser.},
    VOLUME = {181},
     PAGES = {97--119},
 PUBLISHER = {Cambridge Univ. Press, Cambridge},
      YEAR = {1993},
}

@article {LW11,
    AUTHOR = {Lemmens, Bas and Walsh, Cormac},
     TITLE = {Isometries of polyhedral {H}ilbert geometries},
   JOURNAL = {J. Topol. Anal.},
    VOLUME = {3},
      YEAR = {2011},
    NUMBER = {2},
     PAGES = {213--241},
}

@article {MT15,
    AUTHOR = {Matveev, Vladimir S. and Troyanov, Marc},
     TITLE = {Isometries of two dimensional {H}ilbert geometries},
   JOURNAL = {Enseign. Math.},
    VOLUME = {61},
      YEAR = {2015},
    NUMBER = {3-4},
     PAGES = {453--460},
}

@phdthesis {Spe14,
    AUTHOR = {Speer, Timothy},
     TITLE = {Isometries of the {H}ilbert metric},
    SCHOOL = {University of California, Santa Barbara},
      YEAR = {2014},
}

@article {Wal18,
    AUTHOR = {Walsh, Cormac},
     TITLE = {Gauge-reversing maps on cones, and {H}ilbert and {T}hompson isometries},
   JOURNAL = {Geom. Topol.},
    VOLUME = {22},
      YEAR = {2018},
    NUMBER = {1},
     PAGES = {55--104},
}

@article {Mol09,
    AUTHOR = {Moln\'{a}r, Lajos},
     TITLE = {Thompson isometries of the space of invertible positive operators},
   JOURNAL = {Proc. Amer. Math. Soc.},
    VOLUME = {137},
      YEAR = {2009},
    NUMBER = {11},
     PAGES = {3849--3859},
}

@misc {Bos12,
    AUTHOR = {Bosch\'{e}, Aurelien},
     TITLE = {Symmetric cones, the {H}ilbert and {T}hompson metrics},
      NOTE = {Preprint, arXiv:1207.3214},
      YEAR = {2012},
}

@article {LRW19,
    AUTHOR = {Lemmens, Bas and Roelands, Mark and Wortel, Marten},
     TITLE = {Hilbert and {T}hompson isometries on cones in {JB}-algebras},
   JOURNAL = {Math. Z.},
    VOLUME = {292},
      YEAR = {2019},
    NUMBER = {3-4},
     PAGES = {1511--1547},
       DOI = {10.1007/s00209-018-2144-8},
}

@article {RW19,
    AUTHOR = {Roelands, Mark and Wortel, Marten},
     TITLE = {Hilbert isometries and maximal deviation preserving maps on {JB}-algebras},
   JOURNAL = {Adv. Math.},
    VOLUME = {352},
      YEAR = {2019},
     PAGES = {836--861},
       DOI = {10.1016/j.aim.2019.06.027},
}

@book {HOS84,
    AUTHOR = {Hanche-Olsen, Harald and St{\o}rmer, Erling},
     TITLE = {Jordan operator algebras},
    SERIES = {Monographs and Studies in Mathematics},
    VOLUME = {21},
 PUBLISHER = {Pitman (Advanced Publishing Program), Boston, MA},
      YEAR = {1984},
}
\endgroup

\end{document}